\documentclass{article}
\usepackage{fullpage}

\usepackage{amsthm}
\usepackage{amsmath}
\usepackage{amssymb}
\providecommand{\abs}[1]{\left|#1\right|}
\providecommand{\norm}[1]{\lVert#1\rVert}
\providecommand{\re}[1]{\text{Re}(#1)}
\providecommand{\im}[1]{\text{Im}(#1)}
\usepackage{graphicx}
\usepackage{epstopdf}

\usepackage{caption}
\usepackage{subcaption}

\graphicspath{{images/}}
\usepackage{authblk}

\usepackage{mathtools} %for cases
\newtheorem*{remark*}{Remark}
\newtheorem{remark}{Remark}

\newtheorem*{example*}{Example}
\newtheorem{example}{Example}

\newtheorem*{definition*}{Definition}

\newtheorem*{theorem*}{Theorem}

\newtheorem*{proposition*}{Proposition}

\newcounter{saveeqn}
\newcommand{\alpheqn}{\setcounter{saveeqn}{\value{equation}}
\stepcounter{saveeqn}\setcounter{equation}{0}
\renewcommand{\theequation}{\mbox{\arabic{saveeqn}\alph{equation}}}}

\newcommand{\resetalpheqn}{\setcounter{equation}{\value{saveeqn}}
\renewcommand{\theequation}{\arabic{equation}}}

\providecommand{\keywords}[1]{\textbf{\textit{Keywords:}} #1}

\usepackage{xcolor}

\begin{document}
\pagestyle{plain}
\title{The Numerical Unified Transform Method for the Nonlinear Schr\"odinger equation on the half-line}

\author[${\S}$]{Xin Yang \thanks{Corresponding author.}}
\author[$\dagger$]{Bernard Deconinck}
\author[$\ddag$]{Thomas Trogdon}
%\author[$\dagger$]{Xin Yang}
\affil[$$]{Department of Applied Mathematics, University of Washington, Seattle, WA 98195}
\affil[$$]{$^{\dag}$deconinc@uw.edu, $^{\ddag}$trogdon@uw.edu, $^{\S}$yangxin@uw.edu}
\maketitle

\begin{abstract}
We implement the Numerical Unified Transform Method to solve the Nonlinear Schr\"odinger equation on the half-line. For so-called linearizable boundary conditions, the method solves the half-line problems with comparable complexity as the Numerical Inverse Scattering Transform solves whole-line problems. In particular, the method computes the solution at any $x$ and $t$ without spatial discretization or time stepping.  Contour deformations based on the method of nonlinear steepest descent are used so that the method's computational cost does not increase for large $x,t$ and the method is more accurate as $x,t$ increase. Our ideas also apply to some cases where the boundary conditions are not linearizable.
\end{abstract}
\keywords{integrable partial differential equations; Numerical Unified Transform Method; method of nonlinear steepest descent}

\section{Introduction}
\begin{figure}
  \centering
  \begin{subfigure}{0.45\textwidth}
  \includegraphics[width=\textwidth]{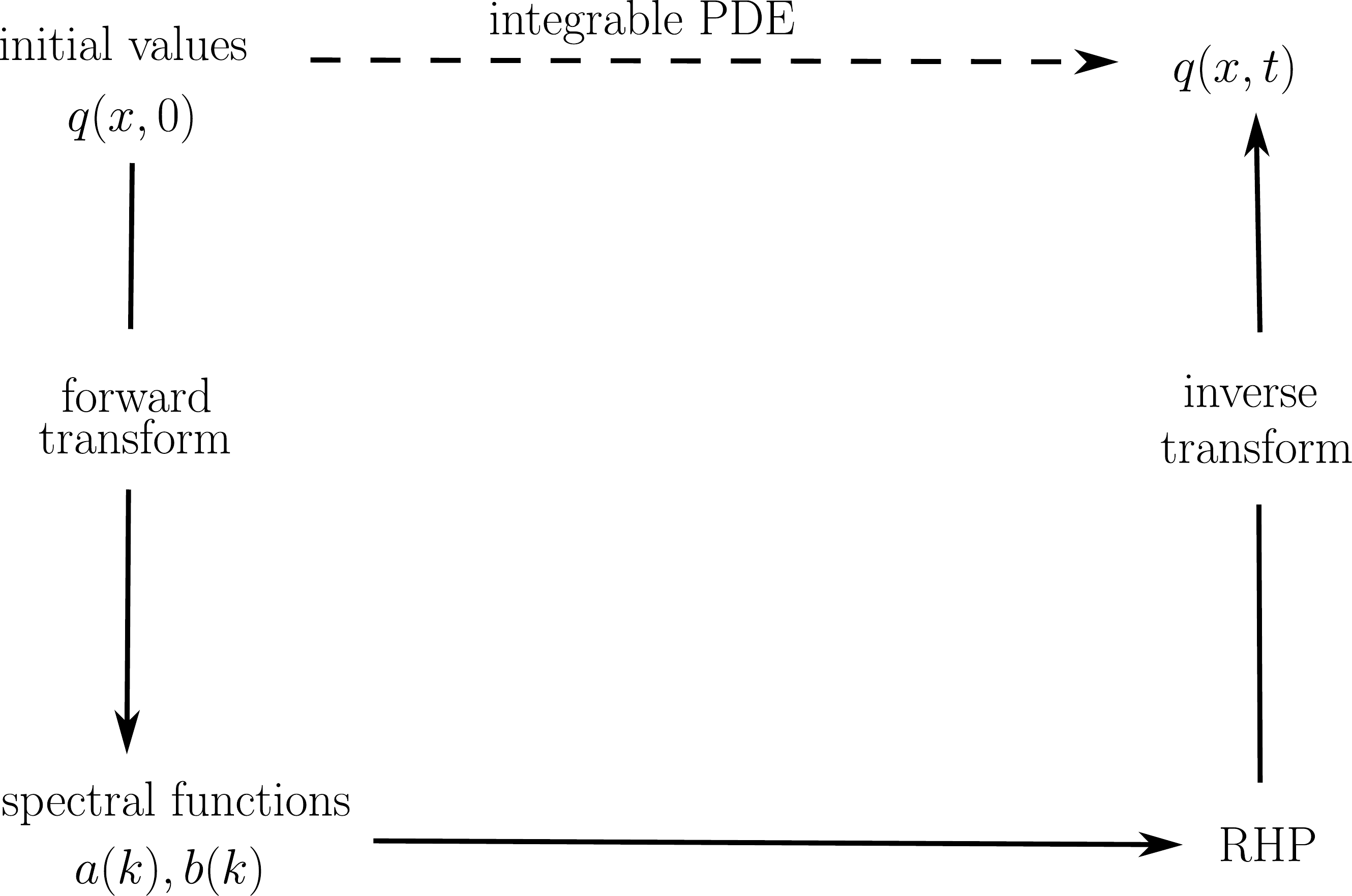}
  \caption{}
  \end{subfigure}
  \quad\quad
  \begin{subfigure}{0.45\textwidth}
  \includegraphics[width=\textwidth]{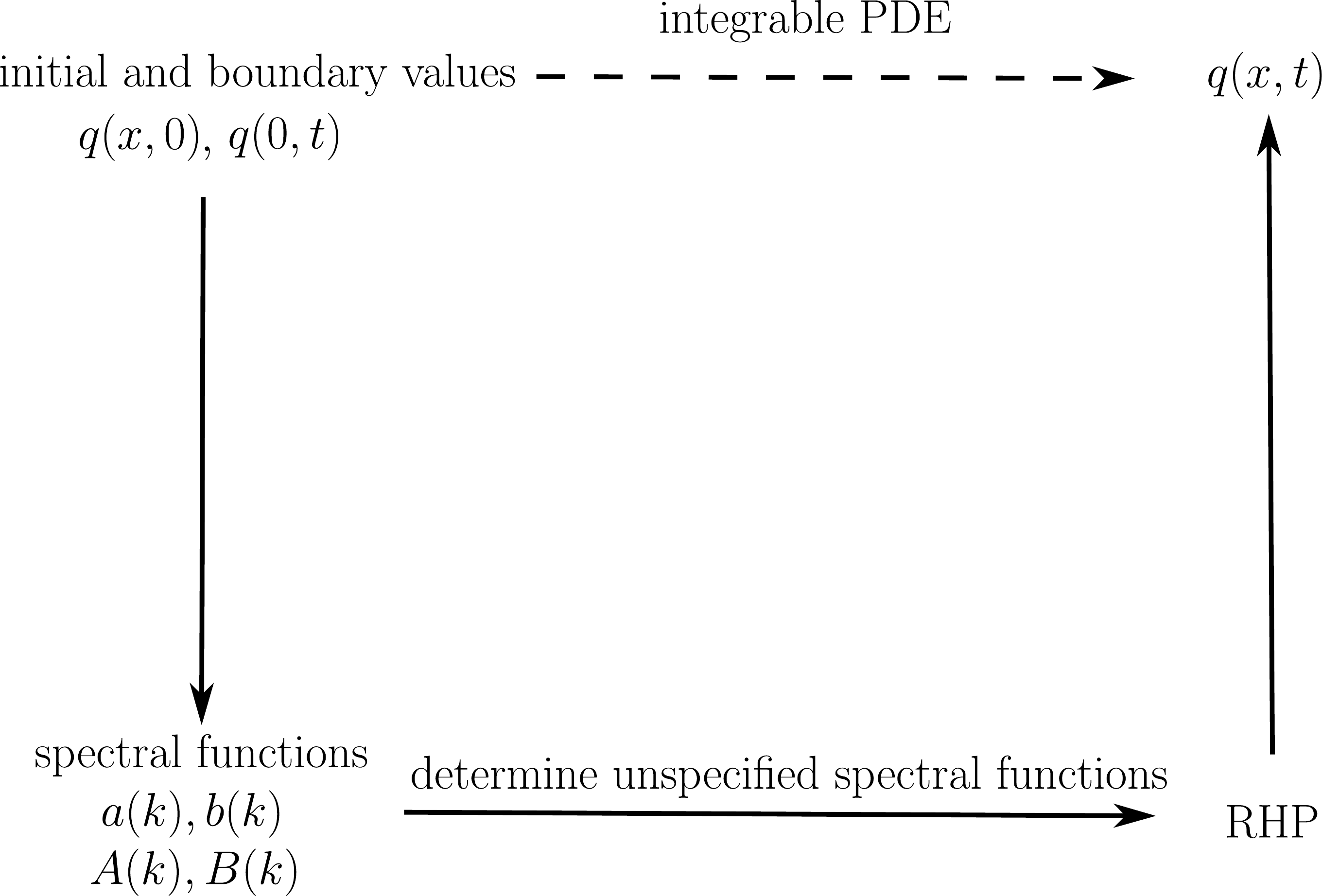}
  \caption{}
  \end{subfigure}
  \caption{A comparison of the schematics of the IST and the UTM. Panel (a) is the diagram of the solution process for the IST for solving an integrable PDE on the whole line. Panel (b) is the diagram of the solution process of the UTM for solving a linear or nonlinear integrable PDE on the half-line. The inverse problems in both the IST and the UTM are formulated as Riemann-Hilbert problems. Dashed lines denote evolution via the integrable PDEs in question. Solid lines denote the steps of the IST and the UTM that can be accomplished by solving a linear problem.}
  \label{diagram}
\end{figure}

In 1997, Fokas developed the Unified Transform Method (UTM) for nonlinear integrable partial differential equations (PDEs) on the half-line \cite{fokas1997}. The UTM is a generalization of the well-known Inverse Scattering Transform (IST) \cite{akns1974} to initial-boundary-value problems (IBVPs). Figure \ref{diagram} illustrates the schematics of the IST and the UTM.
In 2012, Deconinck, Olver and Trogdon developed the Numerical Inverse Scattering Transform (NIST) for the initial-value problems (IVPs) of the Korteweg-de Vries (KdV) equation \cite{trogdon2012a} and the modified Korteweg-de Vries (mKdV) equation \cite{trogdon2012a}. The NIST has since been successfully applied to the IVPs of other integrable systems such as the focusing and defocusing Nonlinear Schr\"odinger (NLS) equations \cite{trogdon2012b}, the Toda lattice \cite{bilman2017} and the sine-Gordon equation \cite{yang2019}.

As a hybrid analytical-numerical method based on the IST, the NIST differs in many aspects from traditional numerical PDE methods. It has the following features:
\begin{enumerate}
\item  The method gives the solution at given $(x,t)$ {\bf without time-stepping or spatial discretization}. \label{feature1}
\item  The method is {\bf spectrally accurate} in the sense that the error at fixed $(x,t)$, $E_{\text{NUTM}}(N,x,t)=\mathcal{O}(1/N^l)$  for any integer $l$, where $N$ is the number of arithmetic operations. \label{feature2}
\item  The method is {\bf uniformly accurate} in the sense that the computational cost to compute the solution at a point $(x,t)$ with given accuracy remains bounded for large $x,t$. \label{feature3}
\item  The method only requires some decay and regularity assumptions on the initial and boundary data. No closed-form expressions for the scattering data are required. \label{feature4}
\item  The method does not artificially truncate infinite physical domains. \label{feature5}
\item  The solution steps require only the solution of linear problems.
\end{enumerate}

Feature 3, which concerns uniform accuracy, comes from the use of the method of nonlinear steepest descent for Riemann-Hilbert Problems (RHPs) \cite{deift1993}. As shown in Figure \ref{diagram}, both the IST and the UTM use RHPs for the inverse transforms. It is natural to ask if it is possible to numerically implement the UTM for IBVPs with the same features. Before we discuss the nonlinear case, it is worth pointing out that although the UTM was originally developed for nonlinear integrable PDEs, it offers a way to analyze IBVPs for linear PDEs, giving the solution in terms of contour integrals in the complex plane~\cite{deconinck2014,fokasbook}. Thus a clear understanding of the differences and difficulties of the UTM for linear problems is necessary for the study of nonlinear problems. For linear evolution equations, there have been numerical methods developed based on the UTM and applied to the heat equation $q_t=q_{xx}$ on the half-line \cite{flyer2008,fokas2009} and on finite intervals \cite{papa2009}, to the Stokes equations $q_t\pm q_{xxx}=0$ on the half-line \cite{flyer2008} and on finite intervals \cite{kesici2018}, and to the advection-diffusion equation $q_t+q_x=q_{xx}$ on the half-line \cite{barros2019}. These applications are implemented with features 1, 2, 5 and 6, but they generally do not possess features 3, 4. Using proper contour deformations and techniques for oscillatory integrals, we recently implemented the Numerical Unified Transform Method (NUTM) for the heat equation, the linear Schr\"odinger (LS) equation and the linear KdV equation with advection on the half-line with features~1-6~\cite{yang2020}.

A major difference in the UTM for nonlinear integrable IBVPs is that the determination of the unknown boundary data, or transforms of this data, from the given initial and boundary conditions is difficult. For general boundary conditions, the unknown boundary data satisfies a system of nonlinear Volterra integral equations \cite{pelloni2015}. We believe that this complication is unavoidable when solving problems with general boundary conditions but we do not consider this further here. In this paper, we show that for linearizable boundary conditions (see Section \ref{sec_linearizablebc} for the definition) the NUTM applies to the NLS equation on the half-line in the same way as for the whole-line problem with all features 1-6. This is indeed expected as these cases can be treated by appropriate spatial reflection. But importantly, in other cases, if the spectral functions are known, the same procedure applies to compute some solutions in the nonlinearizable case.

This paper is organized as follows. Section 2 gives a brief overview of the UTM for the NLS equation on the half-line. In Section 3, we discuss the NUTM applied to the NLS equation with linearizable boundary conditions. In Section 4, we consider the solutions with prescribed spectral functions. The corresponding boundary conditions are not necessarily linearizable. In Section 5, we discuss the asymptotics of the spectral functions to improve accuracy for small $x,t$.

\section{The Unified Transform Method for the NLS equation on the half-line}
In this section we describe the UTM applied to the NLS equation on the half-line. A complete discussion is given in \cite{fokasbook} and \cite{fokas2005}.
\subsection{The Lax pair}
The NLS equation
\begin{align}
iq_t+q_{xx}+2\lambda \abs{q}^2 q =0, ~ \lambda=\pm 1,
\label{nls}
\end{align}
is integrable with the associated Lax pair
\alpheqn
\begin{align}
\mu_x+ik[\sigma_3,\mu]&=Q(x,t)\mu, \label{laxpairx} \\
\mu_t+2ik^2[\sigma_3,\mu]&=\tilde{Q}(x,t,k)\mu, \label{laxpairt}
\end{align}
\resetalpheqn
where $\sigma_3=\mbox{diag}(1,-1)$, $[A,B]=AB-BA$ and
\begin{align*}
Q(x,t)=\left[\begin{matrix}
0 & q(x,t) \\
-\lambda \overline{q(x,t)} & 0
\end{matrix}
\right],~
\tilde{Q}(x,t,k)=2kQ-iQ_x\sigma_3+i\lambda\abs{q}^2 \sigma_3.
\end{align*}
Here $\lambda=\pm 1$ gives the focusing/defocusing NLS equation ($\lambda$ in \cite{fokas2005} is $-\lambda$ here). The compatibility of (2a) and (2b), $\mu_{xt}=\mu_{tx}$, requires $q(x,t)$ satisfying~(\ref{nls}). Using the Lax pair (\ref{laxpairx}),(\ref{laxpairt}) we define
\begin{align}
W(x,t,z)=d(e^{i(kx+2k^2t)\hat{\sigma}_3} \mu(x,t,k))=e^{i(kx+2k^2t)\hat{\sigma}_3}\left(Q\mu(x,t,k) dx+\tilde{Q}\mu(x,t,k) dt\right),
\label{dform}
\end{align}
where
\begin{align*}
\hat{\sigma}_3 A=[\sigma_3, A], ~ e^{\hat{\sigma}_3}A= e^{\sigma_3}A e^{-\sigma_3}.
\end{align*}
Requiring that $W$ is closed implies that $q(x,t)$ satisfies (\ref{nls}). An integral equation for a solution of the Lax pair (\ref{laxpairx}),(\ref{laxpairt}) is obtained by integrating the differential form
\begin{align}
\mu(x,t,k)=I+\int_{(x_*,t_*)}^{(x,t)} e^{-i(kx+2k^2 t) \hat{\sigma}_3} W(\xi,\tau, k ),
\label{musol}
\end{align}
where $I$ is the $2 \times 2$ identity matrix, $x,t\in (0,\infty)$ and $x_*, t_*\in [0,\infty]$. Using particular choices of $(x_*, t_*)$, particular solutions are constructed and used to define the so-called spectral functions.

\subsection{The spectral functions}
Assume that $q(x,t)$ solves (\ref{nls}) for $x,t>0$ with the initial values
\begin{align*}
q(x,0)=q_0(x), ~x \geq 0,
\end{align*}
and the boundary values
\begin{align*}
q(0,t)=g_0(t), ~t \geq 0,\\
q_x(0,t)=g_1(t), ~t \geq 0.
\end{align*}
In general, we do not know both boundary functions $g_0(t)$ and $g_1(t)$. The NLS equation on the half-line is well-posed with either $g_0(t)$ or $g_1(t)$ specified \cite{bu2000,carroll1991,himonas2019}. We further assume that the initial condition $q_0$ is in $\mathcal{S}_{\alpha}(\mathbb{R}^+)$, the Schwartz class functions restricted to the positive half-line with exponential decay rate $\alpha >0 $:
\[
 \mathcal{S}_\alpha(\mathbb{R}^+)=\left\{f \in \mathcal{S}(\mathbb{R})|_{\mathbb{R}^+},~  \exists \alpha'>\alpha>0~\text{:}~\sup_{x \in \mathbb{R}^+} e^{\alpha' x} \abs{f(x)} <\infty  \right\}.
 \]
Following \cite{fokas2005}, in this section, the boundary functions $g_0$ and $g_1$ are assumed to be smooth functions on $[0,T]$ and when $T=\infty$, $g_0$ and $g_1$ are assumed to be in $\mathcal{S}(\mathbb{R})|_{\mathbb{R}^+}$ which is sufficient to define the spectral functions. For the numerical examples we consider, this sufficient condition may not be satisfied. Different assumptions on the boundary values are used for problems discussed in Section 3, 4 and 5 so long as the spectral functions can be computed.

Using the conventions in \cite{fokas2005}, $\mu_1,\mu_2$ and $\mu_3$ are defined using $(x_*,t_*)=(0,T)$, $(x_*,t_*)=(0,0)$ and $(x_*,t_*)=(\infty, t)$, respectively.
For $T<\infty$ the spectral functions $s(k)$ and $S(k,T)$ are defined by
\begin{align*}
s(k)=\mu_3(0,0,k),~  S(k,T)=[e^{2ik^2 T \hat{\sigma}_3} \mu_2 (0,T,k)]^{-1}.
\end{align*}
With this choice, $s(k)$ depends only on the initial values $q_0$ and $S(k,T)$ depends only on the boundary values $g_0$ and $g_1$ on the interval $[0,T]$.
There is an alternate definition $S(k,\infty)=\mu_1(0,0,k)$, which is more convenient when $T=\infty$.

\subsection{Properties of the spectral functions}
Using the symmetries of $Q$ and $\tilde{Q}$, the spectral functions have the form
\begin{align*}
s(k)=\left[\begin{matrix}
\overline{a(\overline{k})} & b(k) \\
-\lambda \overline{b(\overline{k})} & a(k)
\end{matrix}\right],~
S(k,T)=\left[\begin{matrix}
\overline{A(\overline{k},T)} & B(k,T) \\
-\lambda \overline{B(\overline{k},T)} & A(k,T)
\end{matrix}\right].
\end{align*}
Moreover, since $Q$ and $\tilde{Q}$ are traceless,
\begin{align*}
\text{det }s(k)=\text{det }S(k,T)=1,
\end{align*}
which implies that
\begin{align}
a(k)\overline{a(\overline{k})} + \lambda b(k) \overline{b(\overline{k})} &=1, ~ k\in \mathbb{R},\\
A(k,T)\overline{A(\overline{k},T)} + \lambda B(k,T) \overline{B(\overline{k},T)} &=1, ~ k\in \mathbb{C}\,\,(k\in \mathbb{R}\cup i\mathbb{R} ~\text{if } T=\infty). \label{Sdetcondition}
\end{align}
It is convenient to characterize the spectral functions $a(k)=\phi_2(0,k)$ and $b(k)=\phi_1(0,k)$ using linear Volterra integral equations:
\alpheqn
\begin{align}
\phi_1(x,k)&=-\int_x^{\infty} e^{-2ik(x-y)}q_0(y) \phi_2(y,k)dy, \label{phi1}\\
\phi_2(x,k)&=1-\lambda \int_x^{\infty} \bar{q}_0(y) \phi_1(y,k)dy. \label{phi2}
\end{align}
\resetalpheqn
\\
If $q_0\in \mathcal{S}(\mathbb{R})|_{\mathbb{R}^+}$, $a(k)$ and $b(k)$ are analytic for $\im{k}>0 $. With the additional assumption on exponential decay $q_0\in \mathcal{S}_{\alpha}(\mathbb{R}^+)$, $a(k)$ and $b(k)$ are analytic in a larger region that contains $\im{k}\geq -\alpha/2$.

For $T<\infty$, the spectral functions $\overline{A(\overline{k},T)}=\Phi_2(T,k)$ and $-e^{-4ik^2T}B(k,T)=\Phi_1(T,k)$ are defined using a different set of linear Volterra integral equations:
\alpheqn
\begin{align}
\label{phi3}
\Phi_1(t,k)&=\int_0^{t} e^{-4ik^2(t-\tau)} \left( \tilde{Q}_{11} \Phi_1 + \tilde{Q}_{12} \Phi_2   \right) (\tau,k) d\tau, \\
\label{phi4}
\Phi_2(t,k)&=1+\int_0^{t} \left( \tilde{Q}_{21} \Phi_1 + \tilde{Q}_{22} \Phi_2   \right) (\tau,k) d\tau.
\end{align}
\resetalpheqn
Therefore $A(k,T)$ and $B(k,T)$ are entire and bounded in $\im{k^2}\geq  0$.

For $T=\infty$, the spectral functions $A(k,\infty)=\Phi_2(0,k)$ and $B(k,\infty)=\Phi_1(0,k)$ are defined by yet another set of linear Volterra integral equations:
\alpheqn
\begin{align}
\label{phi5}
\tilde{\Phi}_1(t,k)&=-\int_t^{\infty} e^{-4ik^2(t-\tau)} \left( \tilde{Q}_{11} \tilde{\Phi}_1 + \tilde{Q}_{12} \tilde{\Phi}_2   \right) (\tau,k) d\tau, \\
\label{phi6}
\tilde{\Phi}_2(t,k)&=1-\int_t^{\infty} \left( \tilde{Q}_{21} \tilde{\Phi}_1 + \tilde{Q}_{22} \tilde{\Phi}_2   \right) (\tau,k) d\tau.
\end{align}
\resetalpheqn
Therefore $A(k,\infty)$ and $B(k,\infty)$ are analytic for $\im{k^2}>  0$ and bounded for $\im{k^2}\geq  0$.

Since $A(k,T), B(k,T)$ are computed from overdetermined boundary data ($\tilde{Q}$ depends on both $q$ and $q_x$), the spectral functions $a(k)$, $b(k)$, $A(k,T)$, and $B(k,T)$ are not independent. This is clear because the Dirichlet to Neumann map depends on the initial data. The integral of the $1$-form (\ref{dform}) along the boundary of the domain $(x,t)\in (0,\infty)\times (0,T)$ must vanish and we arrive at the global relation connecting the information from the initial value and the boundary values in terms of the spectral data,
\begin{align}
a(k)B(k,T)-b(k)A(k,T)=e^{4ik^2 T} c^{+}(k,T),~ \im{k}\geq 0,
\label{globalrelation}
\end{align}
where $c^+(k,T)$ is:
 \begin{itemize}
 \item an undetermined function analytic for~$\im{k}>0$,
 \item continuous and bounded for $\im{k}\geq 0$, and
 \item $c^+(k,T)=O(1/k)$, as $k\rightarrow \infty$ for $\im{k}>0$.
\end{itemize}
If $T=\infty$, the global relation reduces to
\begin{align}
a(k)B(k,\infty)-b(k)A(k,\infty)=0,~ \im{k}\geq 0,\re{k}\geq 0.
\label{globalrelationinf}
\end{align}

\begin{definition*}[an admissible set of functions \cite{fokas2005}]
Given $q_0 \in \mathcal{S}(\mathbb{R}^+)$, the pair $\{g_0,g_1\}$ of smooth functions on $[0,T]$ or $[0,\infty)$ (if $T=\infty$) is an admissible set of functions with respect to $q_0$ if the following conditions are satisfied:
\begin{enumerate}
\item
    The associated spectral functions $\{a,b,A,B\}$ satisfy the global relation (\ref{globalrelation}) for $T<\infty$ or (\ref{globalrelationinf}) for $T=\infty$.
\item
    The functions $q_0$, $g_0$ and $g_1$ are compatible at $x=t=0$, {\it i.e.}, $g_0(0)=q_0(0)$, $g_1(0)=q'_0(0)$.
    More equation-dependent conditions may be imposed if more regularity of the solution $q(x,t)$ is desired.
\end{enumerate}
\end{definition*}

\begin{remark}
The spectral functions $s(k),{S(k,T)}$ are nonlinear transforms of the corresponding initial and boundary values. For the LS equation (when $\lambda=0$), $\mu_3(x,0,k)=\left(
  \begin{array}{cc}
    \mu_{11} & \mu_{12} \\
    \mu_{21} & \mu_{22} \\
  \end{array}
\right)$ satisfies (\ref{laxpairx}),
\begin{align}
\left(
  \begin{array}{cc}
    \mu_{11} & \mu_{12} \\
    \mu_{21} & \mu_{22} \\
  \end{array}
\right)_x+2ik\left(
  \begin{array}{cc}
    0 & \mu_{12} \\
    -\mu_{21} & 0 \\
  \end{array}
\right)&=\left(
  \begin{array}{cc}
    0 & q(x,0) \\
    0 & 0 \\
  \end{array}
\right)\left(
  \begin{array}{cc}
    \mu_{11} & \mu_{12} \\
    \mu_{21} & \mu_{22} \\
  \end{array}
\right), \label{lslaxpairx}
\end{align}
with $\mu_3(\infty,0,k)=I$. Solving (\ref{lslaxpairx}) yields
\begin{align}
\mu_3(x,0,k)=\left(
  \begin{array}{cc}
    1 & -e^{-2ikx}\int_x^{\infty} e^{2ik \xi} q(\xi,0)d\xi\\
    0 & 1 \\
  \end{array}
\right),
\end{align}
and the spectral function $s(k)$ is given by
\begin{align}
s(k)=\mu_3(0,0,k)=\left(
  \begin{array}{cc}
    1 & -\int_0^{\infty} e^{2ik \xi} q(\xi,0)d\xi\\
    0 & 1 \\
  \end{array}
\right).
\end{align}
Therefore $b(k)=-\int_0^{\infty} e^{2ik \xi} q(\xi,0)d\xi$ is the Fourier transform of the initial condition on the half-line and $a(k)=1$. On the other hand, $\mu_2(0,t,k)=\left(
  \begin{array}{cc}
    \mu_{11} & \mu_{12} \\
    \mu_{21} & \mu_{22} \\
  \end{array}
\right)$ satisfies (\ref{laxpairt}),
\begin{align}
\left(
  \begin{array}{cc}
    \mu_{11} & \mu_{12} \\
    \mu_{21} & \mu_{22} \\
  \end{array}
\right)_t+4ik^2\left(
  \begin{array}{cc}
    0 & \mu_{12} \\
    -\mu_{21} & 0 \\
  \end{array}
\right)&=\left(
  \begin{array}{cc}
    0 & 2kq(0,t)+iq_x(0,t) \\
    0 & 0 \\
  \end{array}
\right)\left(
  \begin{array}{cc}
    \mu_{11} & \mu_{12} \\
    \mu_{21} & \mu_{22} \\
  \end{array}
\right), \label{lslaxpairt}
\end{align}
with $\mu_2(0,0,k)=I$.  Solving (\ref{lslaxpairt}) yields
\begin{align}
\mu_2(0,T,k)=\left(
  \begin{array}{cc}
    1 & e^{-4ik^2 T}\int_0^{T} e^{4ik^2 \eta} \left(2kq(0,\eta)+iq_x(0,\eta)\right) d\eta\\
    0 & 1 \\
  \end{array}
\right),
\end{align}
and the spectral function $S(k,T)$ is given by
\begin{align}
S(k,T)=[e^{2ik^2 T \hat{\sigma}_3} \mu_2 (0,T,k)]^{-1}=
  \left(
  \begin{array}{cc}
    1 & -\int_0^{T} e^{4ik^2 \eta} \left(2kq(0,\eta)+iq_x(0,\eta)\right) d\eta\\
    0 & 1 \\
  \end{array}
\right).
\end{align}
Therefore $B(k,T)=-\int_0^{T} e^{4ik^2 \eta}  \left(2kq(0,\eta)+iq_x(0,\eta)\right) d\eta$ is a sum of Fourier-type transforms of the boundary data on the interval $[0,T]$ and $A(k,T)=1$. The global relation (\ref{globalrelation}) becomes
\[
\int_0^{\infty} e^{2ik \xi} q(\xi,0)d\xi -\left( i \int_0^{T} e^{4ik^2 \eta}  q_x(0,\eta) d\eta + 2k \int_0^{T} e^{4ik^2 \eta}  q(0,\eta) d\eta  \right)=e^{4ik^2 T} c^+(k,T)  ,~ \im{k}\geq 0,
\]
where
\[
c^+(k,T)= \int_0^{\infty} e^{2ik \xi} q(\xi,T)d\xi.
\]
This is exactly the same as the global relation in \cite{deconinck2014} which is obtained using Green's theorem.
\end{remark}

\begin{remark}[The nonlinear Volterra integral equation] In \cite{fokas2012}, it is shown that for the Dirichlet problem, the unknown Neumann value is given by
\begin{align*}
g_1(t)=\frac{2}{\pi i} \int_{\partial D_3} (k \chi_1(t,k)+ig_0(t)) dk +\frac{2g_0(t)}{\pi} \int_{\partial D_3} \chi_2(t,k)dk - \frac{4}{\pi i} \int_{\partial D_3} ke^{-4ik^2 t} \frac{b(-k)}{a(-k)} \overline{\Phi_2 (t,-\overline{k})}dk,
\end{align*}
for $0<t<T$ where $\partial D_3$ is the boundary of the third quadrant in the complex plane with counterclockwise orientation and
\begin{align*}
\chi_j(t,k)=\Phi_j(t,k)-\Phi_j(t,-k), \quad j=1,2, \quad 0<t<T, \quad k\in \mathbb{C}.
\end{align*}
For convenience, suppose that $a(k)$ does not have zeros in the upper half-plane. Plugging the equation of $g_1(t)$ back into (\ref{phi3}) and (\ref{phi4}) yields a nonlinear Volterra integral equation for $\Phi_1(t,k)$ and $\Phi_2(t,k)$ that depend only on known data. However, since $\Phi_j(t,k)$ requires all the values of $\Phi_1(s,k), 0<s<t, k\in \partial D_3$, this is a fundamentally nonlinear problem that depends on two continuous variables. As discussed in Section \ref{sec_overdeterminedbc}, even if one can solve the nonlinear Volterra integral equation for the unknown $g_1(t)$, the slow decay of $g_1(t)$ will affect the overall accuracy when computing the solution to the NLS equation. Perturbative methods for the nonlinear Volterra integral equation have been studied in \cite{fokas2012} and \cite{lenells2012} using related equations.
\end{remark}

\subsection{The Riemann-Hilbert problem}
\begin{figure}
  \centering
  \includegraphics[width=0.5\textwidth]{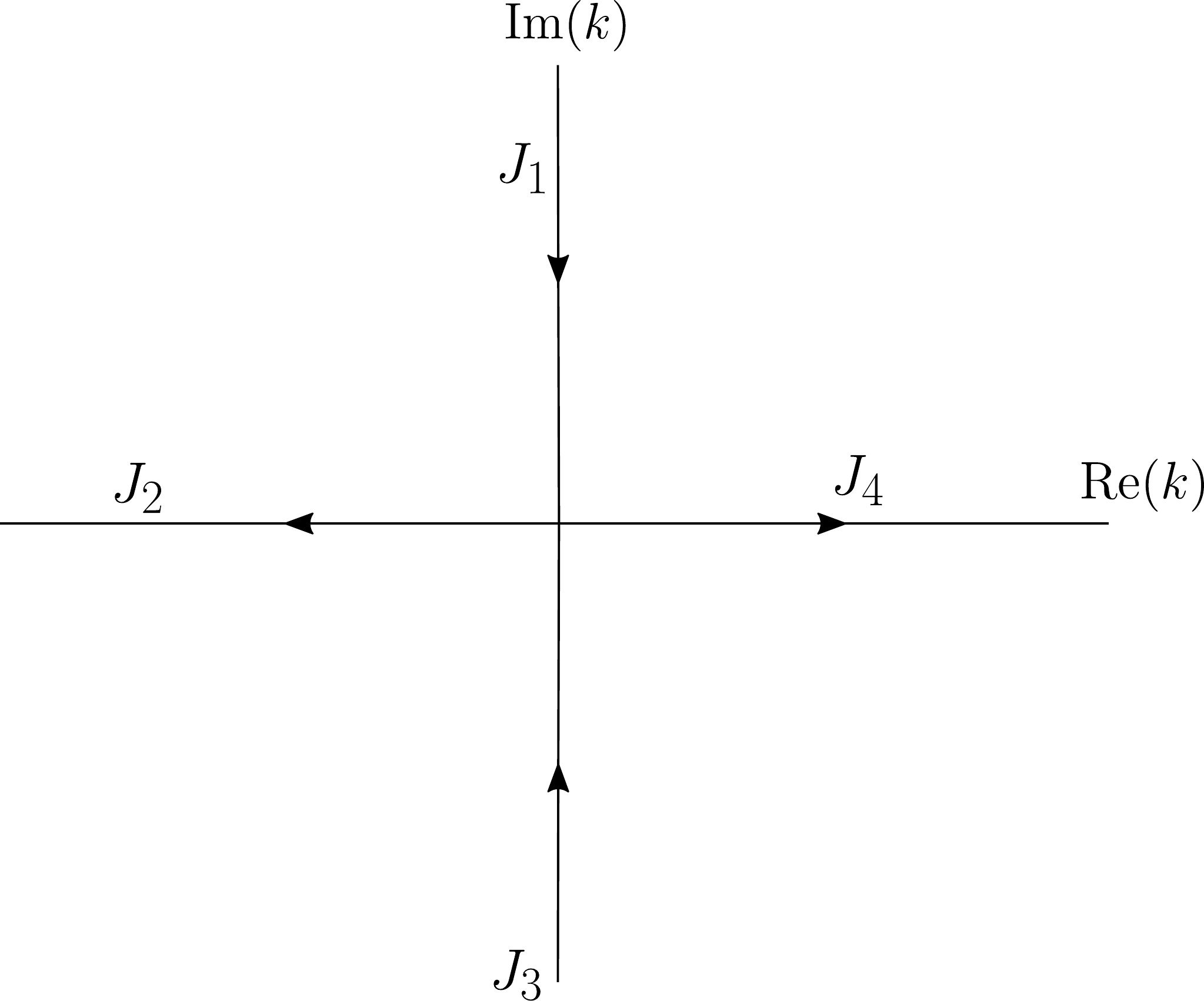}\\
  \caption{The jump contour of RHPs (\ref{rhp1}) and (\ref{rhp2}). The contour consists of the real and imaginary axes.} \label{rhp1contour}
\end{figure}

With the spectral functions $a(k),b(k)$ and $A(k,T),B(k,T)$ defined in the previous sections, we obtain $q(x,t)$ by solving the following RHP.
\begin{theorem*}[\cite{fokas2005}]
Suppose that $q_0$ and an admissible set of functions $\{g_0,g_1\}$ with respect to $q_0$ are given. The spectral functions $a(k),b(k)$ are defined via (\ref{phi1},\ref{phi2}), $A(k,T),B(k,T)$ are defined via (\ref{phi3},\ref{phi4}) for $T<\infty$ and $A(k,\infty),B(k,\infty)$ are defined via (\ref{phi5},\ref{phi6}) for $T=\infty$. Assume that
\begin{itemize}
\item If $\lambda=1$, $a(k)$ has at most $n$ simple zeros $\{p^a_j\}_{j=1}^{n}$ in region $\mathbb{C}^+ \backslash i\mathbb{R}$. Let $n_1$ be the number of zeros in the open first quadrant $\arg p^a_j\in (0,\pi/2), j=1,\dots,n_1,$ and therefore $\arg p^a_j \in (\pi/2,\pi), j=n_1+1,\dots,n$.
\item If $\lambda=1$, for both $T<\infty$ and $T=\infty$, then
\begin{align}
d(k,T)=a(k)\overline{A(\overline{k},T)}+\lambda b(k)\overline{B(\overline{k},T)} &,~ \text{arg } k\in [\pi/2,\pi], \label{ddefinition}
\end{align}
has at most $n_2$ simple zeros in the second quadrant $\{p^d_j\}_{j=1}^{n_2}$, where $\arg p^d_j\in (\pi/2,\pi)$, $j=1,\dots,n_2$.
\end{itemize}
For both $T<\infty$ and $T=\infty$, the following $2\times 2$ matrix RHP for $\Phi(k;x,t)$ has a unique solution:
\begin{enumerate}
\item $\Phi(k;x,t)$ is sectionally meromorphic for $k \in \mathbb{C}\backslash\{\mathbb{R} \cup i\mathbb{R}\}$ (sectionally analytic if $\lambda=-1$).
\item $\Phi(k;x,t)$ has continuous boundary values \footnote{For $k\in \mathcal{C}$ where $\mathcal{C}$ is an oriented contour, define $\Phi_{\pm}(k;x,t)$ to be the limit of $\Phi(k';x,t)$ as $k'\rightarrow k$ nontangentially from the right $(+)$ or the left $(-)$.} on the cross $k\in \mathbb{R}\cup i\mathbb{R}$ with orientation as shown in Figure~\ref{rhp1contour}, and
\begin{align}
\Phi_{+}(k;x,t)=\Phi_{-}(k;x,t)J(k;x,t).
\label{rhp1}
\end{align}
The jump matrix $J(k;x,t)$ is given by
\begin{align}
J(k;x,t)= \begin{cases} J_4(k;x,t), & \mbox{arg } k = 0, \\
                        J_1(k;x,t), & \mbox{arg } k = \frac{\pi}{2},\\
                        J_2(k;x,t), & \mbox{arg } k = \pi,\\
                        J_3(k;x,t), & \mbox{arg } k = \frac{3\pi}{2},
    \end{cases}
    \label{rhp1jump}
\end{align}
where
\begin{align*}
J_1(k;x,t)=&\left[\begin{matrix}
1 & 0 \\
-\Gamma(k,T)e^{2i\theta(k;x,t)} & 1
\end{matrix}\right],\\
J_3(k;x,t)=&\left[\begin{matrix}
1 & -\lambda \overline{\Gamma(\overline{k},T)}e^{-2i\theta(k;x,t)} \\
0 & 1
\end{matrix}\right],\\
J_4(k;x,t)=&\left[\begin{matrix}
1+\lambda\gamma(k)\overline{\gamma(\overline{k})}  & \gamma(k)e^{-2i\theta(k;x,t)} \\
\lambda\overline{\gamma(\overline{k})}e^{2i\theta(k;x,t)} & 1
\end{matrix}\right],\\
J_2(k;x,t)=&J_1(k;x,t) J_4^{-1}(k;x,t) J_3 (k;x,t)\\
=& \left[\begin{matrix}
1  & -\left(\lambda \overline{\Gamma(\overline{k},T)} + \gamma(k)\right)e^{-2i\theta(k;x,t)} \\
-\left(\lambda \overline{\gamma(\overline{k})} + \Gamma(k,T)\right)e^{2i\theta(k;x,t)} & 1 +\left(\lambda \overline{\Gamma(\overline{k},T)} + \gamma(k)\right)\left(\lambda \overline{\gamma(\overline{k})} + \Gamma(k,T)\right)
\end{matrix}\right],
\end{align*}
and
\begin{align}
\theta(k;x,t)=kx+2k^2 t &, ~ k\in \mathbb{C}, \label{rhp1theta} \\
\gamma(k)=\frac{b(k)}{\overline{a(\overline{k})}} &, ~  k\in \mathbb{R},  \label{rhp1lowergamma}\\
\Gamma(k,T)=-\frac{\lambda \overline{B(\overline{k},T)}}{a(k)d(k,T)} &, ~ \text{arg }  k\in \left[\frac{\pi}{2},\pi\right]. \label{rhp1uppergamma}
\end{align}

\item If $\lambda=1$, the first column of $\Phi(k;x,t)$ has simple poles at $\{p_j^a\}_{j=1}^{n_1}$ and $\{p_j^d\}_{j=1}^{n_2}$. The second column of $\Phi(k;x,t)$ has simple poles at $\{\overline{p_j^a}\}_{j=1}^{n_1}$ and $\{\overline{p_j^d}\}_{j=1}^{n_2}$. The associated residues satisfy the relations:
\alpheqn
\begin{align}
\label{res1}
\text{Res}_{k=p^a_j}\,\Phi(k;x,t)&=\lim_{k\rightarrow p^a_j} \Phi(k;x,t) \left(
                                                            \begin{array}{cc}
                                                              0 & 0 \\
                                                              \frac{1}{a'(p^a_j)b(p^a_j)} e^{2i\theta(p^a_j;x,t)} & 0 \\
                                                            \end{array}
                                                          \right), \quad  j=1,\dots, n_1, \\
\label{res2}
\text{Res}_{k=\overline{p^a_j}}\,\Phi(k;x,t)&=\lim_{k\rightarrow \overline{p^a_j}} \Phi(k;x,t) \left(
                                                            \begin{array}{cc}
                                                              0 & \frac{-\lambda}{\overline{a'(p^a_j)b(p^a_j)}} e^{-2i\theta(\overline{p^a_j};x,t)} \\
                                                              0 & 0 \\
                                                            \end{array}
                                                          \right), \quad  j=1,\dots, n_1,\\
\label{res3}
\text{Res}_{k=p^d_j}\,\Phi(k;x,t)&=\lim_{k\rightarrow p_j^d} \Phi(k;x,t) \left(
                                                            \begin{array}{cc}
                                                              0 & 0 \\
                                                              -\frac{\lambda \overline{B(\overline{p^d_j}})}{a(p^d_j)d'(p^d_j)} e^{2i\theta(p^d_j;x,t)} & 0 \\
                                                            \end{array}
                                                          \right), \quad  j=1,\dots, n_2, \\
\label{res4}
\text{Res}_{k=\overline{p^d_j}}\,\Phi(k;x,t)&=\lim_{k\rightarrow \overline{p^d_j}} \Phi(k;x,t) \left(
                                                            \begin{array}{cc}
                                                              0 & \frac{B(\overline{p^d_j})}{\overline{a(p^d_j)d'(p^d_j)}} e^{-2i\theta(\overline{p^d_j};x,t)} \\
                                                              0 & 0 \\
                                                            \end{array}
                                                          \right), \quad  j=1,\dots, n_2.
\end{align}
\resetalpheqn
\item $\Phi(k;x,t)=I+O(1/k)$ as $k\rightarrow \infty$.
\end{enumerate}
Then (as is shown in \cite{fokas2005})
\[
q(x,t)=2i \lim_{k\rightarrow \infty} \left(k \Phi(k;x,t) \right)_{12},\]
solves the NLS equation with
\[
q_x(x,t)=\lim_{k\rightarrow \infty} \left\{4  \left(k^2 \Phi(k;x,t)\right)_{12} + 2i q(x,t) \left(k  \Phi(k;x,t)\right)_{22} \right\},
\]
\[
q(x,0)=q_0(x),\quad q(0,t)=g_0(t) \quad \text{and} \quad q_x(0,t)=g_0(t).
\]
\end{theorem*}

\begin{figure}
  \centering
  \includegraphics[width=0.9\textwidth]{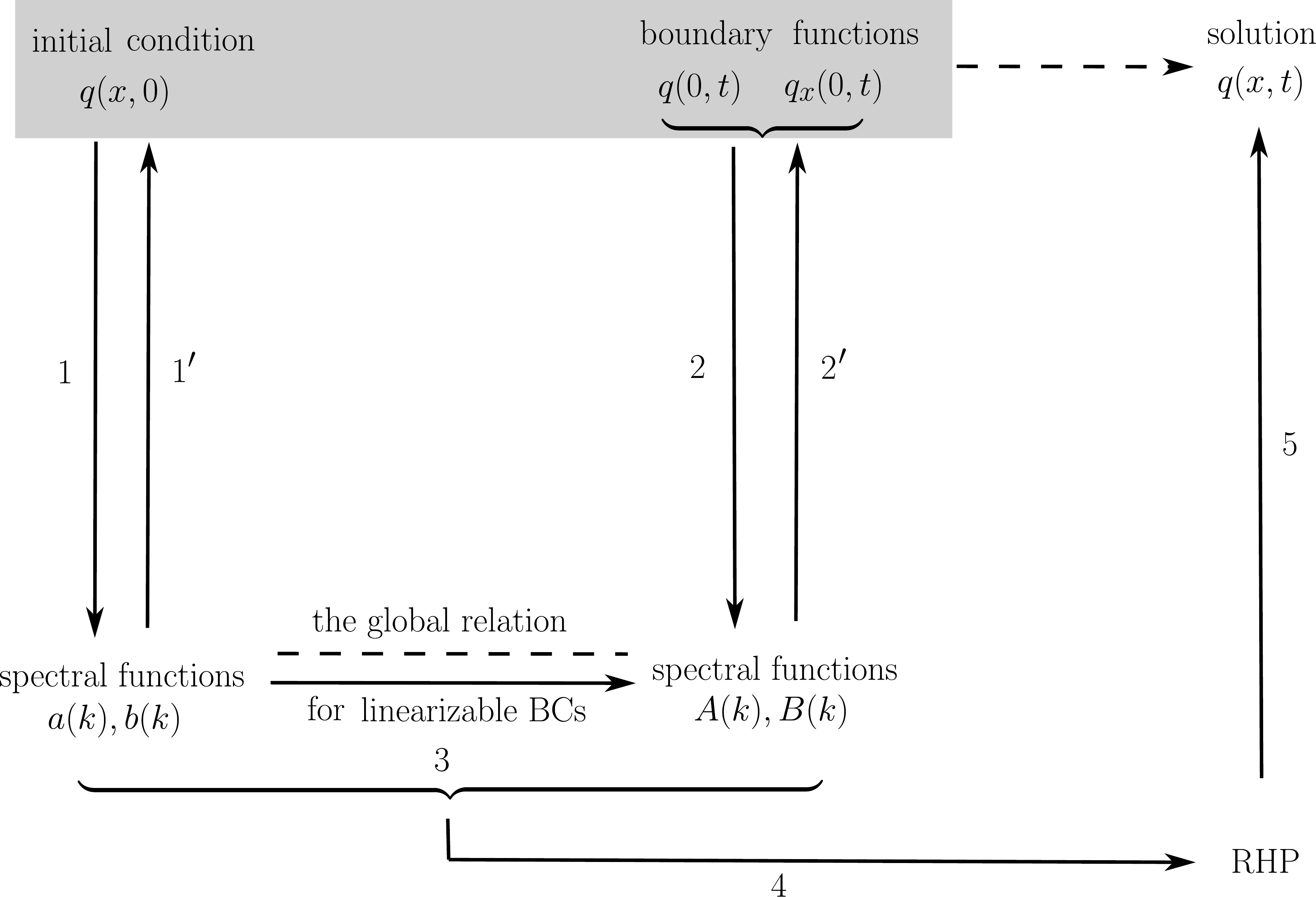}
  \caption{A diagram describing the use of the UTM to solve integrable PDEs on the half-line. The path along 1, 4, 5 is the same as the IST for IVPs. In this paper, we focus on three paths: (i) the path along 1-5 for problems with linearizable boundary conditions, (ii) the path along 1, 2, 4, 5 for problem with overdetermined but compatible initial and boundary values, and (iii) the path along 4, 5 for problems with specified spectral functions.}
  \label{diagram2}
\end{figure}

Figure \ref{diagram2} shows a detailed diagram of the UTM applied to nonlinear integrable PDEs. Steps $1$ and~$1'$ denote the invertible transforms between the initial condition $q(x,0)=q_0(x)$ and the spectral functions $\{a(k),b(k)\}$. Step $2$ and $2'$ denote the invertible transforms between the pair of boundary functions $\{q(0,t)=g_0(t)$, $q_x(0,t)=g_1(t)\}$ and the spectral functions $\{A(k,T),B(k,T)\}$ for $T<\infty$ or $T=\infty$. Step~$3$ denotes obtaining $\{A(k,T),B(k,T)\}$ within a special class of boundary conditions known as linearizable boundary conditions as discussed in Section \ref{sec_linearizablebc}. Step~$4$ denotes the construction of the associated RHP using the spectral functions incorporating time dependence. Step $5$ denotes the inverse transform to get the solution $q(x,t)$ by solving the RHP. In the following sections, we show examples of the solutions that can be efficiently computed with different types of given data. In Section~\ref{sec_linearizablebc}, we give an example with linearizable boundary conditions following steps $1$-$5$. In Section 4.1, we show an example with an admissible set of functions $g_0,g_1 \in \mathcal{S}_\beta(\mathbb{R}^+)$ with respect to $q_0 \in \mathcal{S}_\alpha(\mathbb{R}^+)$. This follows steps $1,2,4,5$ in Figure \ref{diagram2}. In Section 4.2, we show an example with specified spectral functions that follows steps $4,5$ in Figure \ref{diagram2}, which can also be understood as an application of the dressing method to construct solutions to the NLS equation \cite{fokasbook}.

\begin{remark}
The defocusing NLS equation on the half-line does not have soliton solutions \cite{lenells2016}. More precisely, $a(k)$ does not have zeros for $\im{k}\geq 0$ and $d(k,T)$ does not have zeros for arg$(k)\in [\pi/2,\pi]$.
\end{remark}

\begin{remark}
When $t=0$, the RHP (\ref{rhp1}) reduces to a RHP that depends only on $a(k)$ and $b(k)$ by deforming $J_1$ and $J_3$ to the negative real line. The global relation is not needed in the deformation. When $x=0$, the RHP (\ref{rhp1}) reduces to a RHP that depends only on $A(k,T)$ and $B(k,T)$ on the cross $k\in \mathbb{R}\cup i\mathbb{R}$ but the reduction requires the use of the global relation as well as (\ref{Sdetcondition}).
\end{remark}

\subsection{Algorithms for computing the spectral data}
\label{computespectral}
The goal of the numerical computation for the forward transform is to obtain: (i) the evaluation of the spectral functions along the (deformed) jump contour, (ii) in the focusing case, the zeros of $a(k)$ and $d(k,T)$ as well as the related residues. \\
\\
{\bf (i) Continuous spectral data.}
For convenience, we describe the algorithm for the computation of the spectral data in the case $T=\infty$. Assuming $q_0$ and an admissible set of functions $g_0$ and $g_1$ are given, we compute the spectral functions $\{a(k),b(k)\}$ using the differential equation form of (\ref{phi1},\ref{phi2}). We compute $A(k,\infty),B(k,\infty)$ using the differential equation form of (\ref{phi5},\ref{phi6}). %The second column in both sets of (\ref{laxpairx}) and (\ref{laxpairt})
Both sets of equations are in the form
\begin{align*}
\left(
  \begin{array}{c}
    y_1(s,k) \\
    y_2(s,k) \\
  \end{array}
\right)_s + \mathcal{M}_1(k) \left(
  \begin{array}{c}
    y_1(s,k) \\
    y_2(s,k)+1 \\
  \end{array}
\right)
=
 \mathcal{M}_2(s,k) \left(
  \begin{array}{c}
    y_1(s,k) \\
    y_2(s,k)+1 \\
  \end{array}
\right).
\end{align*}
They are solved by a Chebyshev collocation method \cite{battles2004} on $[0,L]$ with vanishing boundary condition at $y_1(L,k)=y_2(L,k)=0$ for sufficiently large $L$. A detailed discuss of the Chebyshev collocation method solving this type of the equations can be found in \cite{trogdon2012a} and \cite{tombook}.
\\
\\
{\bf (ii) Discrete spectral data.} In \cite{lenells2016}, it is shown that the zeros of $a(k)$ in the upper half-plane are the $L^2(\mathbb{R},\mathbb{C}^{2\times 2})$ eigenvalues of the operator
    \[
    \mathcal{L}=i\sigma_3\partial_x-i\sigma_3 Q_e,
    \]
    where
    \[
    q_e(x)= \begin{cases}
q(x,0), & x\geq 0, \\
0, & x<0, \\
\end{cases}
    \quad\text{and}\quad
    Q_e=\left[\begin{matrix}
0 & q_e \\
-\lambda \overline{q_e} & 0
\end{matrix}
\right].
\]
The eigenvalues are obtained using the Floquet-Fourier-Hill method \cite{deconinck2006}. Though the Floquet-Fourier-Hill method does not achieve spectral accuracy due to the possible discontinuity of the potential $q_e(x)$ at $x=0$, it provides initial guesses for Newton's method. The residue conditions require evaluating $c^a_j=\frac{1}{a'(p^a_j)b(p^a_j)}$ at the zeros $\{p_j^a\}_1^{n_1}$ where $a'(p^a_j)$ is computed using Cauchy's integral formula.

Similarly, it is shown in \cite{lenells2016} that the zeros of $d(k,T)$ satisfy the same eigenvalue problem except that
\[
    q_e(x)= \begin{cases}
q(x,T), & x\geq 0, \\
0, & x<0. \\
\end{cases}\]
The potential $q_e$ depends on the unknown solution $q(x,T)$, and pure root-finding algorithms are needed to find zeros of $d(k,T)$.
Once the zeros $\{p^d_j\}_1^{n_2}$ are obtained, $c^d_j=-\frac{\lambda \overline{B(\overline{p^d_j}})}{a(p^d_j)d'(p^d_j)}$ are computed using Cauchy's integral formula.\\

\section{The NLS equation with linearizable boundary conditions}
\label{sec_linearizablebc}
\subsection{Linearizable boundary conditions}
Obtaining $A(k,T)$ and $B(k,T)$ is non-trivial since they are defined in terms of the Dirichlet {\em and} Neumann data, both of which cannot be arbitrarily specified, for a well-posed problem~\cite{bu2000,carroll1991,himonas2019,holmer2005}. In the special case of linearizable boundary conditions, $A(k,T)$ and $B(k,T)$ can be obtained by solving algebraic equations involving $a(k)$ and $b(k)$ without solving (\ref{laxpairt}) (recall that (\ref{laxpairt}) requires boundary functions $g_0(t)$ and $g_1(t)$). The idea is to use the global relation (\ref{globalrelation}) to find extra identities using the symmetries of the dispersion relation similar to how the UTM is applied to linear PDEs \cite{deconinck2014}. For the LS equation, the dispersion relation is $\omega(k)=ik^2$, which is invariant under the mapping $k\rightarrow -k$. Similarly, for the NLS equation, we want to determine a relation between $A(k,T),B(k,T)$ and $A(-k,T),B(-k,T)$. Recall that $A(k,T)$ and $B(k,T)$ are defined in terms of $e^{2ik^2T \hat{\sigma}_3}\mu_2(0,T,k)$. Let $\Phi(t,k)=\mu_2(0,t,k)e^{-2ik^2t\sigma_3}$, then $\Phi(t,k)$ satisfies
\begin{align*}
\Phi_t+2ik^2 \sigma_3 \Phi=\tilde{Q}(t,k) \Phi, ~\Phi(0,k)=I.
\end{align*}
Suppose there exists a $t$-independent, nonsingular matrix $N(k)$ such that
\begin{align}
\Phi(t,-k)=N(k)\Phi(t,k)N^{-1}(k).
\label{lbc1}
\end{align}
More explicitly, (\ref{lbc1}) is equivalent to
\begin{align}
(2ik^2\sigma_3-\tilde{Q}(t,-k))N(k)=N(k)(2ik^2\sigma_3-\tilde{Q}(t,k)).
\label{lbc2}
\end{align}
A necessary condition for the existence of $N(k)$ is that the determinant of $2ik^2\sigma_3-\tilde{Q}(t,k)$ is even in $k$. This implies
\begin{align*}
q(0,t)\overline{q}_x(0,t)-\overline{q}(0,t)q_x(0,t)=0.
\end{align*}
If this condition is satisfied, (\ref{lbc2}) becomes
\begin{align*}
(2kq-iq_x)N_{21}&=\lambda (2k\overline{q}-i\overline{q}_x)N_{12},\\
(2kq+iq_x)N_{11}+(2kq-iq_x)N_{22}&=-2(2ik^2-i\lambda \abs{q}^2)N_{12}.
\end{align*}
In particular, for the homogeneous Robin boundary condition with a real parameter $\rho>0$ (this choice of sign is discussed below in Remark \ref{remark_rbbc})
\begin{align}
q_x(0,t)-\rho q(0,t)=0,
\label{rbc}
\end{align}
we choose $N_{12}=N_{21}=0$ and $(2k-i\rho)N_{22}+(2k+i\rho)N_{11}=0$ so that
\begin{align}
A(k,T)=A(-k,T),~B(k,T)=-\frac{2k+i\rho}{2k-i\rho}B(-k,T), ~k\in \mathbb{C}.
\label{robin}
\end{align}
The results for the homogeneous Dirichlet boundary condition and the homogeneous Neumann boundary condition are obtained by taking $\rho\rightarrow \infty$ and $\rho \rightarrow 0$ respectively.
For instance, for the homogeneous Dirichlet boundary condition, (\ref{robin}) becomes
\begin{align}
A(k,T)=A(-k,T),~B(k,T)=B(-k,T), ~k\in \mathbb{C}.
\label{dirichlet}
\end{align}
With this, we can solve for $A(k)$ and $B(k)$ in terms of $a(k)$ and $b(k)$ explicitly. Indeed, if $T=\infty$, the global relation (\ref{globalrelationinf}) gives the following equation, valid in the first quadrant,
 \begin{align}
 a(k)B(k,\infty)-b(k)A(k,\infty)=0, ~\text{arg } k \in \left[0,\frac{\pi}{2}\right].
\label{eq1}
 \end{align}
Letting $k\rightarrow -k$ in the expression for $\overline{d(\overline{k},\infty)}$ in (\ref{ddefinition}) and using (\ref{dirichlet}), we find a second equation, also valid in the first quadrant:
\begin{align}
A(k,\infty) \overline{a(-\overline{k})} +\lambda B(k,\infty) \overline{b(-\overline{k})}=\overline{d(-\overline{k},\infty)}, ~\text{arg } k \in \left[0,\frac{\pi}{2}\right].
\label{eq2}
\end{align}
Solving both (\ref{eq1}) and (\ref{eq2}) for $A(k,\infty)$ and $B(k,\infty)$ yields
\begin{align}
A(k,\infty)=\frac{a(k) \overline{d(-\overline{k},\infty)}}{\Delta_0(k)}, ~ B(k,\infty)=\frac{b(k) \overline{d(-\overline{k},\infty)}}{\Delta_0(k)}, ~\text{arg } k \in \left[0,\frac{\pi}{2}\right],
\label{AB}
\end{align}
where
\begin{align}
\Delta_0(k)=a(k)  \overline{a(-\overline{k})} +\lambda b(k) \overline{b(-\overline{k})}, ~\text{arg } k \in [0,\pi] .
\label{delta0}
\end{align}
There is no need to solve for $d(k,\infty)$ in terms of $a(k)$ and $b(k)$ since the jump condition in (\ref{rhp1}) depends only on $B(k,\infty)$, $A(k,\infty)$ through the ratio $B(k,\infty)/A(k,\infty)$.
\\
\\
\begin{example}[The RHP associated with homogenous Dirichlet boundary conditions]
With a homogenous Dirichlet boundary condition, we obtain a RHP involving only $a(k)$ and $b(k)$, which are determined solely by the initial condition. We seek a $2\times2$ matrix-valued function $\Phi(k;x,t)$ that satisfies
\begin{align}
\Phi_{+}(k;x,t)=\Phi_{-}(k;x,t)J(k;x,t),
\label{rhp2}
\end{align}
with the jump functions defined in (\ref{rhp1jump}) on the cross $k\in \mathbb{R}\cup i\mathbb{R}$, shown in Figure \ref{rhp1contour}. The only difference is that (\ref{rhp1uppergamma}) becomes
\begin{align}
\Gamma(k)=\frac{-\lambda \overline{b(-\overline{k})}}{a(k) \Delta_0(k)} &, ~\text{arg } k \in [0,\pi],
\end{align}
where $\Delta_0(k)$ is given by (\ref{delta0}).
\end{example}
\begin{remark} In general, $\rho$ in the homogeneous Robin boundary condition (\ref{rbc}) can take any real value \cite{bu2000}. As shown in \cite{its2013}, when $\rho<0$, generically, there are zeros of $a(k)$ on the positive imaginary axis. This requires modifications of our assumptions on the RHPs since we do not allow for poles on the jump contour. Meanwhile, the long time behavior of such solutions at $x=0$ is dominated by oscillatory standing solitons leading to non-decaying boundary data. The choice of the sign of $\rho$ is also related to the possibility of extending the half-line solution to a bounded whole-line solution, see \cite{biondini2012,fokas2004} for further details.
\label{remark_rbbc}
\end{remark}

\subsection{Deformation of the contour based on the method of nonlinear steepest descent}
We use the numerical approach developed in \cite{olver2012,tombook} to solve the RHP (\ref{rhp2}). Uniform accuracy can be obtained using appropriate deformations of the jump contours. Then, the deformed RHP is solved using the Mathematica package {\it RHPackage} developed by Olver \cite{olver2012} with spectral accuracy. The deformations are derived in a similar fashion as the deformations used for the solution of the RHP for the NLS equation on the whole line \cite{trogdon2012b}. The idea is to deform the contour near the saddle point to the steepest descent direction so that the oscillations from the exponential factor $e^{2i\theta(k;x,t)}$ change to exponential decay. The saddle point $k_0$ of the phase $\theta(k;x,t)$ is determined by
\[
\frac{d\theta(k;x,t)}{dk}\Big|_{k=k_0}=0 \Rightarrow k_0= -\frac{x}{4t}.
\]
We write the exponent as
\[
2i\theta(k;x,t)=-\frac{ix^2}{4t}+4it(k-k_0)^2.
%+\mathcal{O}(k-k_0)^3
\]
Thus $e^{2i\theta(k;x,t)}$ is exponentially decaying if $k$ follows a path with arg$(k-k_0)=\pi/4,5\pi/4$. In addition, the deformation of contours requires that the functions $\gamma(k),\Gamma(k)$ are analytic in the neighbourhood of $k_0$. Since $q_0 \in \mathcal{S}_\alpha(\mathbb{R}^+)$, $a(k)$ and $b(k)$ are analytic and bounded for $\im{k}\geq -\alpha/2$. To ensure that the residue condition is outside the region in which the contour is deformed, if $0\leq \min(\{\im{p_j^a}\}_{j=1}^n)\leq \alpha/2$, we redefine $\alpha=\min(\{\im{p_j^a}\}_{j=1}^n)/4$. Therefore, $\gamma(k)$ is bounded and analytic in a strip centered around the real axis with $-\alpha/2 \leq \im{k}\leq \alpha/2$, while $\Gamma(k)$ is bounded and analytic for $\im{k}\geq -\alpha/2$. In some cases, if $\tau(k)$ in (\ref{tau}) vanishes in the strip, then we need to further shrink the width of the strip. Since $\tau(k)\geq 1$ on the real axis, we can always find a valid choice for $\alpha>0$. We introduce the following deformation steps for the contour of the RHP (\ref{rhp2}).\\
\subsubsection{Step 1: deformations based on steepest descent directions}
\begin{figure}
  \centering
  \includegraphics[width=0.8\textwidth]{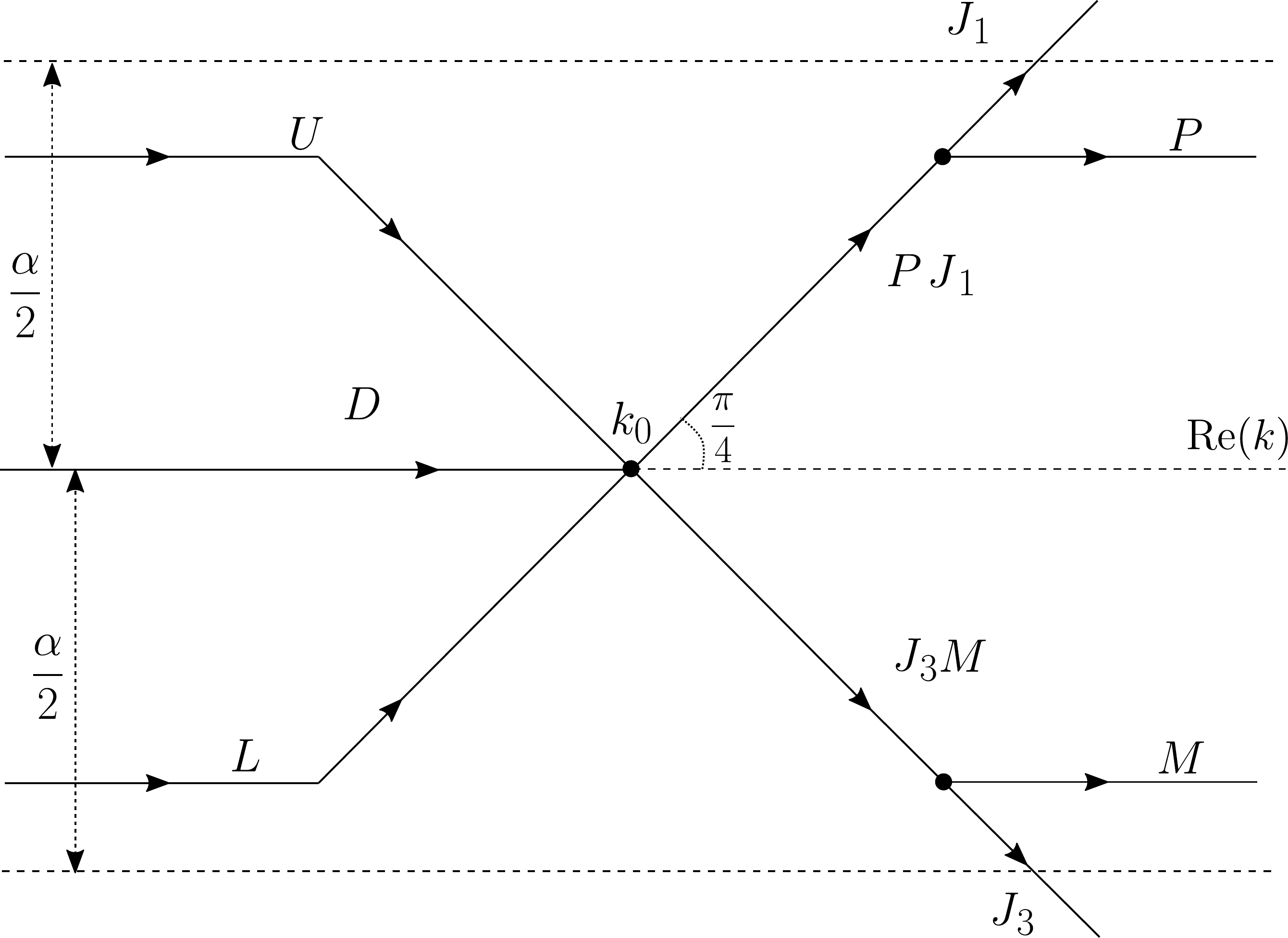}\\
  \caption{The deformed contour for the RHP (\ref{rhp3}) in the complex $k$-plane near the saddle point $k_0$ for the method of nonlinear steepest descent. }
  \label{contour0}
\end{figure}

Let $\mathcal{R}_{k,\theta}=\{k+r e^{i\theta}: r\geq 0\}$. The jump matrix $J_4(k;x,t)$ has the factorization
\[
J_4(k;x,t)=\left[\begin{matrix}
1+\lambda\gamma(k)\overline{\gamma(\overline{k})}  & \gamma(k)e^{-2i\theta(k;x,t)} \\
\lambda\overline{\gamma(\overline{k})}e^{2i\theta(k;x,t)} & 1
\end{matrix}\right]=\left[\begin{matrix}
1  & \gamma(k)e^{-2i\theta(k;x,t)} \\
0 & 1
\end{matrix}\right]\left[\begin{matrix}
1  & 0 \\
\lambda\overline{\gamma(\overline{k})}e^{2i\theta(k;x,t)} & 1
\end{matrix}\right]=MP.
\]
This factorization provides $J_4(k;x,t)$ with decay away from $k_0$ as, by replacing the contour on the real line with two oblique rays starting from $k_0$. Then $M$ approaches the identity matrix exponentially fast in the lower half-plane along $\mathcal{R}_{k_0,7\pi/4}$ and $P$ approaches the identity matrix exponentially fast in the upper half-plane along $\mathcal{R}_{k_0,\pi/4}$. However, for $\re{k}<\re{k_0}$, the exponentials in $M$ and $P$ are growing for $M$ along $\mathcal{R}_{k_0,5\pi/4}$ in the lower half-plane and $P$ along $\mathcal{R}_{k_0,3\pi/4}$ in the upper half-plane. Alternatively, $J_2(k;x,t)$ has the factorization
\begin{align*}
J_2(k;x,t)&=\left[\begin{matrix}
1  & -\left(\lambda \overline{\Gamma(\overline{k})} + \gamma(k)\right)e^{-2i\theta(k;x,t)} \\
-\left(\lambda \overline{\gamma(\overline{k})} + \Gamma(k)\right)e^{2i\theta(k;x,t)} & 1 +\left(\lambda \overline{\Gamma(\overline{k})} + \gamma(k)\right)\left(\lambda \overline{\gamma(\overline{k})} + \Gamma(k)\right)
\end{matrix}\right]\\
&=\left[\begin{matrix}
1  & 0 \\
-\left(\lambda \overline{\gamma(\overline{k})} + \Gamma(k)\right) \frac{ e^{2i\theta(k;x,t)} } {\tau(k)} & 1
\end{matrix}\right]
\left[\begin{matrix}
\frac{1}{\tau(k)}  & 0 \\
0 & \tau(k)
\end{matrix}\right]
\left[\begin{matrix}
1  & -\left(\lambda \overline{\Gamma(\overline{k})} + \gamma(k)\right) \frac{e^{-2i\theta(k;x,t)}}{\tau(k)} \\
0 & 1
\end{matrix}\right]=LDU,
\end{align*}
where
\begin{align}
\tau(k)=
1 + \left(\lambda \overline{\gamma(\overline{k})} + \Gamma(k)\right)
  \left(\lambda \overline{\Gamma(\overline{k})} + \gamma(k)\right) ,~ \im{k}\leq \frac{\alpha}{2}.
\label{tau}
\end{align}
This factorization provides $J_2(k;x,t)$ with decay for increasing $t$ as, for $\re{k}<\re{k_0}$, $L$ approaches the identity matrix exponentially fast along $\mathcal{R}_{k_0,5\pi/4}$ in the lower half-plane and $U$ approaches the identity matrix exponentially fast along $\mathcal{R}_{k_0,3\pi/4}$ in the upper half-plane. We obtain the RHP
\begin{align}
\Phi_{+}(k;x,t)=\Phi_{-}(k;x,t)H(k;x,t),
\label{rhp3}
\end{align}
with jump functions
\begin{align}
H(k;x,t)= \begin{cases}
                        M(k;x,t), & \{k\in \mathcal{R}_{k_0,7\pi/4}: \abs{\im{k}}\leq \alpha/2\}\cup \{k\in \mathbb{C}:k=k_0+\frac{\alpha}{2}e^{-i\frac{\pi}{4}}+s,~ s\geq 0 \}, \\
                        P(k;x,t), & \{k\in \mathcal{R}_{k_0,\pi/4}: \abs{\im{k}}\leq \alpha/2\}\cup \{k\in \mathbb{C}:k=k_0+\frac{\alpha}{2}e^{i\frac{\pi}{4}}+s,~ s\geq 0 \}, \\
                        J_1(k;x,t), & \{k\in \mathcal{R}_{k_0,\pi/4}\},\\
                        L(k;x,t), & \{k\in \mathcal{R}_{k_0,5\pi/4}: \abs{\im{k}}\leq \alpha/2\}\cup \{k\in \mathbb{C}:k=k_0+\frac{\alpha}{2}e^{i\frac{5\pi}{4}}-s,~ s\geq 0 \}, \\
                        D(k;x,t), & \{k\in \mathbb{C}:k=k_0-s,~ s\geq 0 \}, \\
                        U(k;x,t), & \{k\in \mathcal{R}_{k_0,3\pi/4}: \abs{\im{k}}\leq \alpha/2\}\cup \{k\in \mathbb{C}:k=k_0+\frac{\alpha}{2}e^{i\frac{3\pi}{4}}-s,~ s\geq 0 \}, \\
                        J_3(k;x,t), & \{k\in \mathcal{R}_{k_0,7\pi/4}\},\\
    \end{cases}
\end{align}
and the deformed contour, with orientation, is shown in Figure \ref{contour0}.

\subsubsection{Step 2: deformations for uniform accuracy}
\begin{figure}
  \centering
  \includegraphics[width=0.8\textwidth]{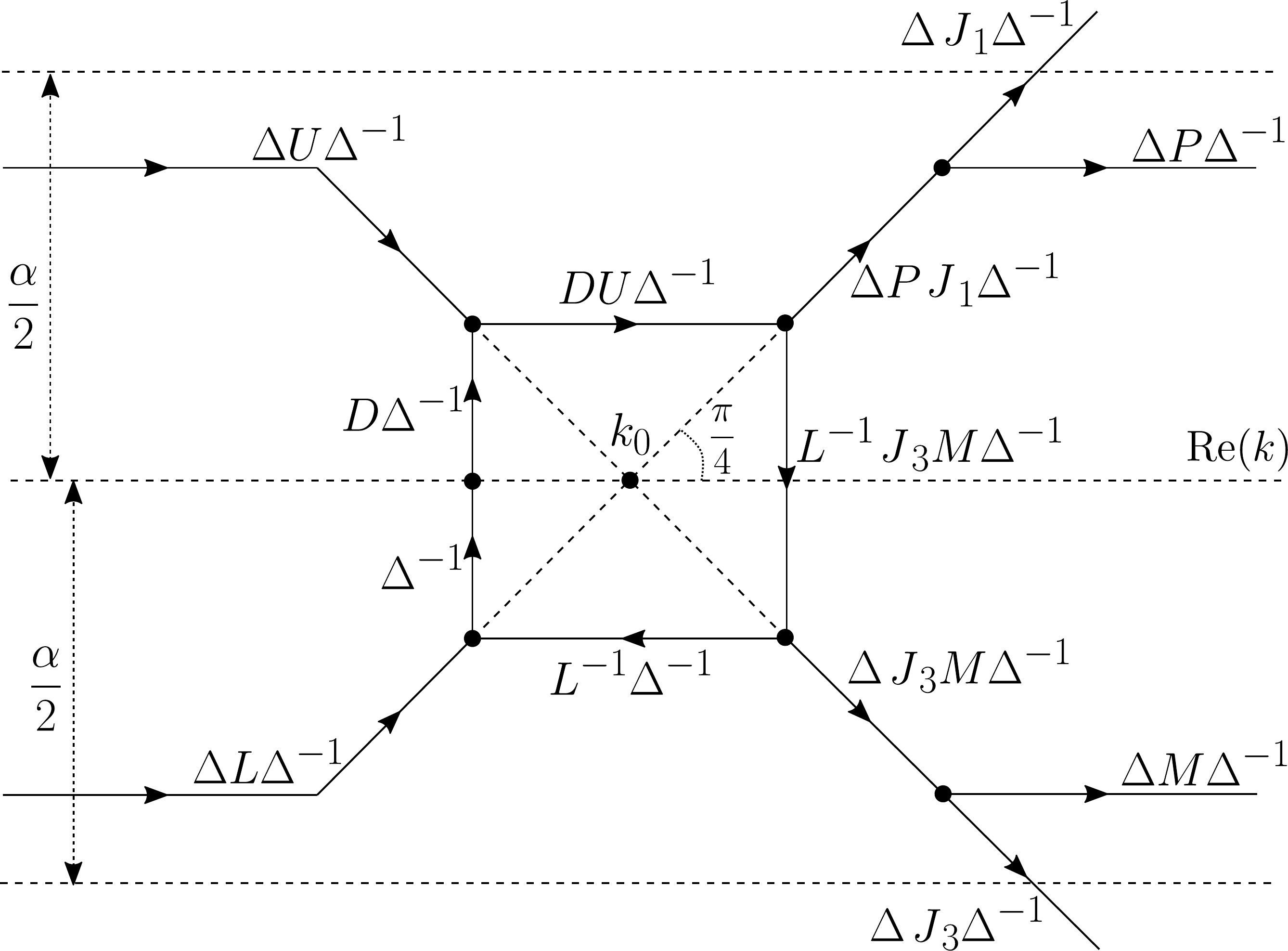}\\
  \caption{The deformed contour for the RHP (\ref{rhp4}) in the complex $k$-plane near the saddle point $k_0$ after removing the jump on the negative real axis. All jumps away from $k_0$ approach the identity exponentially fast as $t\rightarrow \infty$. The length of the side of the square is on the order of $O(1/\sqrt{t})$ for large $t$.}
  \label{contour1}
\end{figure}

Similar to the RHP for the whole-line problem in \cite{trogdon2012b}, the errors for computing the solution of RHP (\ref{rhp3}) are not uniformly small for large time since not all jumps decay to the identity matrix away from the saddle point $k_0$. For large $t$, although the jump matrix $D$ along the negative real axis does not contain oscillatory exponentials, the solution of the RHP (\ref{rhp3}) has increasing oscillations along the jump contour for $k<k_0$ as $t$ grows. Therefore we remove the jump matrix $D$ using conjugation \cite{trogdon2012b}. We introduce the matrix-valued function $\Delta(k,k_0)$,
\begin{align*}
\Delta(k,k_0)= \left[\begin{matrix}
\delta(k,k_0)  & 0 \\
0 & \delta^{-1}(k,k_0)
\end{matrix}\right],
\end{align*}
where
\begin{align*}
\delta(k,k_0)=\exp\left( \frac{1}{2\pi i}\int_{-\infty}^{k_0} \frac{\log \tau(z)}{z-k} dz\right).
\end{align*}
Then $\Psi(k;x,t)=\Phi(k;x,t)\Delta^{-1}(k,k_0) $ is continuous across the real axis for $\re{k}\neq \re{k_0}$, and satisfies
\begin{align*}
\Psi_{+}(k;x,t)=\Psi_{-}(k;x,t) \Delta(k,k_0) H(k;x,t) \Delta^{-1}(k,k_0)=\Psi_{-}(k;x,t) \tilde{H}(k;x,t).
\end{align*}
Since $\delta(k,k_0)$ is singular at $k_0$, lensing is used to avoid the singularity by introducing new jump conditions on a square around $k_0$, as shown in Figure \ref{contour1}. The length of the side of the square is $O(1/\sqrt{t})$ for large $t$. See \cite{trogdon2012b} for further details of the scaling. Summarizing all deformations, we have the RHP
\begin{align}
\Psi_{+}(k;x,t)=\Psi_{-}(k;x,t) \tilde{H}(k;x,t),
\label{rhp4}
\end{align}
where the jump contours are shown in Figure \ref{contour1}.
%$\tilde{J}(k;x,t)$ become
%\begin{align}
%\tilde{J}(k;x,t)= \begin{cases}
%                        \Delta(k,k_0)M(k;x,t) \Delta^{-1}(k,k_0) & \{k\in \mathbb{C}:k=k_0+se^{-i\frac{\pi}{4}},~\frac{\delta}{2} \geq s\geq 0 \}\cup \{k\in \mathbb{C}:k=k_0+\frac{\delta}{2}e^{-i\frac{\pi}{4}}+s,~ s\geq 0 \}, \\
%                        \Delta(k,k_0)P(k;x,t) \Delta^{-1}(k,k_0) & \{k\in \mathbb{C}:k=k_0+se^{i\frac{\pi}{4}},~\frac{\delta}{2} \geq s\geq 0 \}\cup \{k\in \mathbb{C}:k=k_0+\frac{\delta}{2}e^{i\frac{\pi}{4}}+s,~ s\geq 0 \}, \\
%                        \Delta(k,k_0)J_1(k;x,t) \Delta^{-1}(k,k_0) & \{k\in \mathbb{C}:k=k_0+se^{i\frac{\pi}{4}},~s\geq 0 \},\\
%                        \Delta(k,k_0)L(k;x,t) \Delta^{-1}(k,k_0) & \{k\in \mathbb{C}:k=k_0+se^{i\frac{5\pi}{4}},~\frac{\delta}{2} \geq s\geq 0 \}\cup \{k\in \mathbb{C}:k=k_0+\frac{\delta}{2}e^{i\frac{5\pi}{4}}-s,~ s\geq 0 \}, \\
%                        %D(k;x,t) & \{k\in \mathbb{C}:k=k_0-s,~ s\geq 0 \}, \\
%                        \Delta(k,k_0)U(k;x,t) \Delta^{-1}(k,k_0) & \{k\in \mathbb{C}:k=k_0+se^{i\frac{3\pi}{4}},~\frac{\delta}{2} \geq s\geq 0 \}\cup \{k\in \mathbb{C}:k=k_0+\frac{\delta}{2}e^{i\frac{3\pi}{4}}-s,~ s\geq 0 \}, \\
%                        \Delta(k,k_0)J_3(k;x,t) \Delta^{-1}(k,k_0) & \{k\in \mathbb{C}:k=k_0+se^{-i\frac{\pi}{4}},~s\geq 0 \},\\
%    \end{cases}
%\end{align}
%The jump contour of $D$ in the RHP (\ref{rhp3}) is therefore removed.

\subsubsection{Step 3: Adding residue conditions}
For the focusing NLS equation, there is an additional step for the residue conditions (\ref{res1})-(\ref{res4}). By introducing small circles centered at the singularities, and modifying the unknown function $\Psi$ inside the circle, the residue conditions are replaced with jump conditions on the circles \cite{trogdon2012b}. Let $\{z_j\}_{1=j}^{n_1+n_2}$ be the union of the zeros of $a(k)$ and $d(k,T)$ defined by  $z_j=p_j^{a}$, for $1\leq j\leq n_1$ and $z_j=p_{j-n_1}^{d}$, for $n_1+1 \leq j\leq n_1+n_2$. Let $\{c_j\}_{1=j}^{n_1+n_2}$ be defined by  $c_j=c_j^{a}, 1\leq j\leq n_1$ and $c_j=c_{j-n_1}^{d},n_1+1 \leq j\leq n_1+n_2$.
For a residue condition at $k=z_j$ in the upper half-plane,
\begin{align}
  \Psi^+(k;x,t)= \Psi^-(k;x,t) \left(
                                                            \begin{array}{cc}
                                                              1 & 0 \\
                                                              c_j e^{2i\theta(z_j;x,t)}/(k-z_j) & 1 \\
                                                            \end{array}
                                                          \right),
\end{align}
is the jump condition on a circle centered at $z_j$ with radius $\varepsilon$ oriented counterclockwise. The circles need to avoid intersections with contours already present in the RHP. The corresponding residue condition at $k=\overline{z_j}$ in the lower half-plane becomes
\begin{align}
  \Psi^+(k;x,t)= \Psi^-(k;x,t) \left(
                                                            \begin{array}{cc}
                                                              1 & -\overline{c_j} e^{-2i\theta(\overline{z_j};x,t)}/(k-\overline{z_j}) \\
                                                              0 & 1 \\
                                                            \end{array}
                                                          \right),
\end{align}
on a circle centered at $\overline{z_j}$ with radius $\varepsilon$.
Since $\abs{c_j e^{2i\theta(z_j;x,t)}}$ may be unbounded for large $x,t$, we invert this factor through a deformation when $\abs{c_j e^{2i\theta(z_j;x,t)}}>1$. We define the matrix-valued function $\hat{\Psi}(k;x,t)$ by
\begin{align}
\hat{\Psi}(k;x,t)=
\begin{cases}
    \Psi(k;x,t)\left(
        \begin{array}{cc}
          1 & -(k-z_j)/(c_j e^{2i\theta(z_j;x,t)}) \\
          (c_j e^{2i\theta(z_j;x,t)})/(k-z_j) & 0 \\
        \end{array}
      \right)V(k),  &  \text{if } \abs{k-z_j}<\varepsilon, \\
      \Psi(k;x,t)\left(
        \begin{array}{cc}
          0 & (-\overline{c_j} e^{-2i\theta(\overline{z_j};x,t)})/(k-\overline{z_j}) \\
          (k-\overline{z_j})/(\overline{c_j} e^{-2i\theta(\overline{z_j};x,t)}) & 1\\
        \end{array}
      \right)V(k),  &  \text{if } \abs{k-\overline{z_j}}<\varepsilon,\\
      \Psi(k;x,t)V(k),  &  \text{otherwise,}
\end{cases}
\end{align}
where
\[
v(z)=\prod_{j\in K_{x,t}} \frac{k-z_j}{k-\overline{z_j}} \text{\quad and \quad} V(z)=\left(
                    \begin{array}{cc}
                      v(z) & 0 \\
                      0 & 1/v(z) \\
                    \end{array}
                  \right),
\]
for each $j$ in the set $K_{x,t}=\{j: \abs{c_j e^{2i\theta(z_j;x,t)}}>1\}$ containing the indexes of the zeros of $a(k)$ and $d(k,T)$ whose jump matrices need to be inverted.
Then $\hat{\Psi}(k;x,t)$ satisfies the jump conditions
\[
\hat{\Psi}_{+}(k;x,t)=\hat{\Psi}_{-}(k;x,t) V^{-1}(k)\tilde{H}(k;x,t)V(k),
\]
on the same contours as (\ref{rhp4}). In addition, $\hat{\Psi}(k;x,t)$ satisfies the jump conditions circles around $\{z_j\}_{1=j}^{n_1+n_2}$,
\begin{align*}
\hat{\Psi}^{+}(k;x,t)=
\begin{cases}
    \hat{\Psi}^{-}(k;x,t)V^{-1}(k)\left(
        \begin{array}{cc}
          1 & -(k-z_j)/(c_j e^{2i\theta(z_j;x,t)}) \\
          0 & 1 \\
        \end{array}
      \right)V(k),  &  \text{if } \abs{k-z_j}=\varepsilon, \\
     \hat{\Psi}^{-}(k;x,t) V^{-1}(k)\left(
        \begin{array}{cc}
           1& 0 \\
          (k-\overline{z_j})/(\overline{c_j} e^{-2i\theta(\overline{z_j};x,t)}) & 1\\
        \end{array}
      \right)V(k),  &  \text{if } \abs{k-\overline{z_j}}=\varepsilon.
\end{cases}
\end{align*}
With all deformations, the RHP (\ref{rhp4}) is solved using {\it RHPackage}~\cite{olver2012} after truncating the contours along which the jump matrices are close to the identity matrix. In practice, this tolerance is set to $10^{-9}$ unless otherwise stated. We use the same tolerance $10^{-9}$ when solving for the spectral functions. For convenience, we also truncate the contours if they are outside a disk centered at the origin with radius $50$. In most cases, the truncation errors are on the same order of the tolerance since the jump matrix approaches the identity matrix exponentially fast. When $x$, $t$ are small, the truncation error dominates. We discuss how to control the truncation error in Section \ref{sec_fasterdecay}.

\subsection{Numerical results}
\begin{figure}
  \centering
  \includegraphics[width=0.8\textwidth]{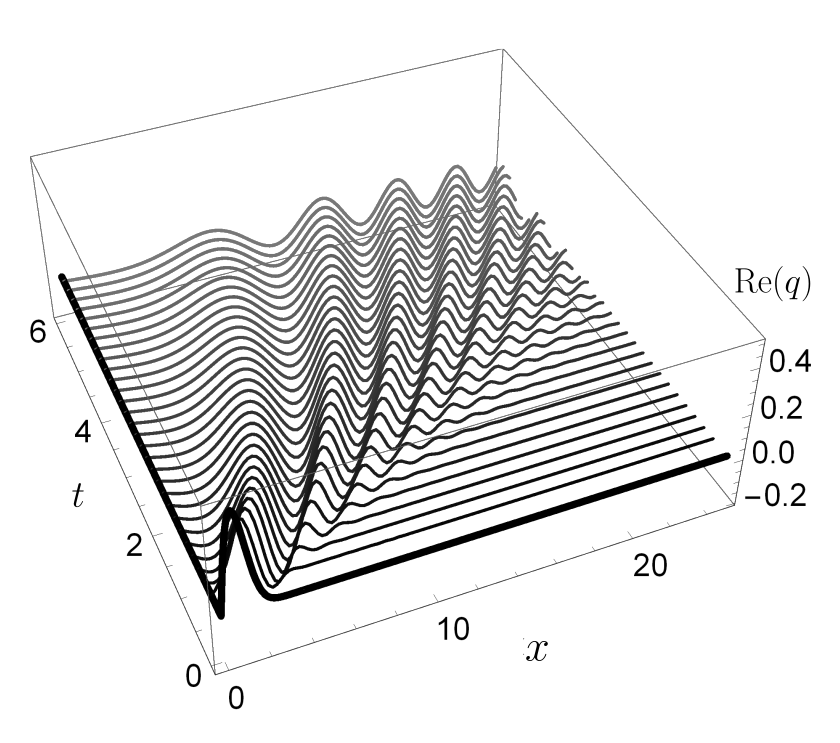}\\
  \caption{The real part of $q(x,t)$ with homogenous Dirichlet boundary condition $q(0,t)=0$ and the initial condition $q(x,0)=x e^{-x^2}$. The thick curves show the initial and boundary conditions.}\label{homodiri}
\end{figure}

In Figure \ref{homodiri}, we plot the real part of the solution to the defocusing NLS equation on the half-line with homogenous Dirichlet boundary condition at $x=0$ and the initial condition $q(x,0)=x e^{-x^2}$. We observe dispersive waves propagating to the right from the localized initial condition. We plot the solution for  $0\leq x\leq 10$ and $0.1\leq t\leq 3.5$. The domain is chosen to be bounded away from $t=0$ as the NUTM is less efficient there. Indeed, when the UTM is applied to linear PDEs on the half-line, the solution formula requires principle-value integrals for its evaluation at $x=t=0$ \cite{yang2020}. This issue occurs with the NLS equation at $x=t=0$ as well. The UTM is well-defined for any $x>0$ or $t>0$ but the numerics suffer from slow convergence when $x,t$ are small, see Section 5.

\section{The NLS equation with nonlinearizable boundary conditions}
\subsection{Overdetermined boundary conditions}
\label{sec_overdeterminedbc}
A possible way to avoid computing unknown boundary conditions is to specify both boundary functions $g_0(t)$ and $g_1(t)$, provided they are admissible with respect to the given initial condition $q_0(x)$. However, there are obstacles to computing the associated solution efficiently:
\begin{enumerate}
\item For a generic whole-line solitonless solution $q(x,t)$ with a nontrivial reflection coefficient $\rho(k)=b(k)/a(k)$, one has $q(0,t)\sim t^{-1/2}$ and $q_x(0,t) \sim t^{-3/2}$ as $t\rightarrow \infty$ \cite{fokas2005}. Therefore, in general, for half-line problems the Dirichlet and Neumann data do not both decay exponentially. This affects the regions where the contours can be deformed. For instance, in the $T=\infty$ case, $A(k,\infty)$ and $B(k,\infty)$ are only guaranteed to be analytic in the first and third quadrants. Therefore, jump contours , which depend on $\tau(k)$ defined in (\ref{tau}) cannot be deformed away from the real axis and the method of nonlinear steepest descent cannot be applied directly. These undeformed contours become highly oscillatory as $t$ increases. For linear PDEs, numerical methods such as Levin's method can be used to compute the oscillatory integrals with high accuracy \cite{yang2020}. For nonlinear integrable PDEs, efficient numerical methods for oscillatory singular integral equations from the RHP are not as well developed \cite{trogdon2015}. A complete discussion of this is beyond the scope of this paper.
\item For pure whole-line soliton solutions that have non-zero velocity, $q(0,t)$ and $q_x(0,t)$ decay exponentially. In this case, the focusing NLS equation allows right-going soliton solutions whose parameters correspond to zeros of $d(k,T)$ in the second quadrant. As discussed in Section \ref{computespectral}, this step requires root-finding algorithms.
\end{enumerate}
There are solutions, with compatible $q_0,g_0,g_1$, that do not suffer from (1) and (2). We can compute these solutions efficiently. Such solutions include left-going singular solutions of the defocusing NLS equation, known as positons, or left-going soliton solutions of the focusing NLS equation.

Such solutions have analytical expressions and are used to demonstrate the accuracy of the NUTM. In fact, unlike for the whole line problem, the jump function (\ref{rhp1jump}) in the half-line problem is non-trivial even when the solution does not contain dispersion. For instance, suppose the initial and boundary values are prescribed by the one-positon solution of the defocusing NLS equation \cite{tombook}
\[
q(x,t)=2\eta e^{-4it(\xi^2-\eta^2)-2ix \xi }\text{ csch}(2\eta(4t\xi +x-x_0)),
\]
where $\xi,\eta,x_0$ are constants. The positon is left-going if $\xi>0$. If $x_0<0$ the singularity is outside the domain for all $t\geq 0$, therefore $q(x,t)$ is exponentially localized and smooth for $x,t\geq 0$. This positon solution corresponds to a simple zero of $a(k)$ at $k_1=\xi+i \eta \tanh(2\eta x_0)$ \cite{kalimeris2010}. The assumptions $\xi>0$ and $x_0<0$ imply that $k_1$ is not in the first quadrant so no residue conditions are required for formulating the associated RHP. Similarly, suppose the initial and boundary values are obtained from the one-soliton solution of the focusing NLS equation \cite{tombook},
\begin{align}
q(x,t)=2\eta e^{-4it(\xi^2-\eta^2)-2ix \xi}\text{ sech}(2\eta(4t\xi +x-x_0)).
\label{onesoliton}
\end{align}
This soliton corresponds to a simple zero of $a(k)$ at $k_1=\xi+i \eta \tanh(2\eta x_0)$. If the initial position of the center of the soliton lies outside of the domain ({\it i.e.,} $x_0<0$) and it is left-going ({\it i.e.,} $\xi<0$), the soliton is represented by the continuous part of the spectral data and no residue conditions are used. On the other hand, for the focusing NLS equation, it is possible to allow $x_0>0$. In this case, $k_1$ is in the first quadrant and residue conditions are required:
\begin{align*}
\text{Res}_{k=k_1}\,\Phi(k;x,t)&=\lim_{k\rightarrow k_1} \Phi(k;x,t) \left(
                                                            \begin{array}{cc}
                                                              0 & 0 \\
                                                              c_1 e^{2i\theta(k_1;x,t)} & 0 \\
                                                            \end{array}
                                                          \right),\\
\text{Res}_{k=\overline{k_1}}\,\Phi(k;x,t)&=\lim_{k\rightarrow \overline{k_1}} \Phi(k;x,t) \left(
                                                            \begin{array}{cc}
                                                              0 & -\overline{c_1} e^{-2i\theta(\overline{k_1};x,t)} \\
                                                              0 & 0 \\
                                                            \end{array}
                                                          \right),
\end{align*}
where $c_1=1/(a'(k_1)b(k_1))$. These considerations can be generalized to $n$-positon and $n$-soliton solutions.

\begin{figure}
  \centering
  \begin{subfigure}{0.45\textwidth}
  \includegraphics[width=\textwidth]{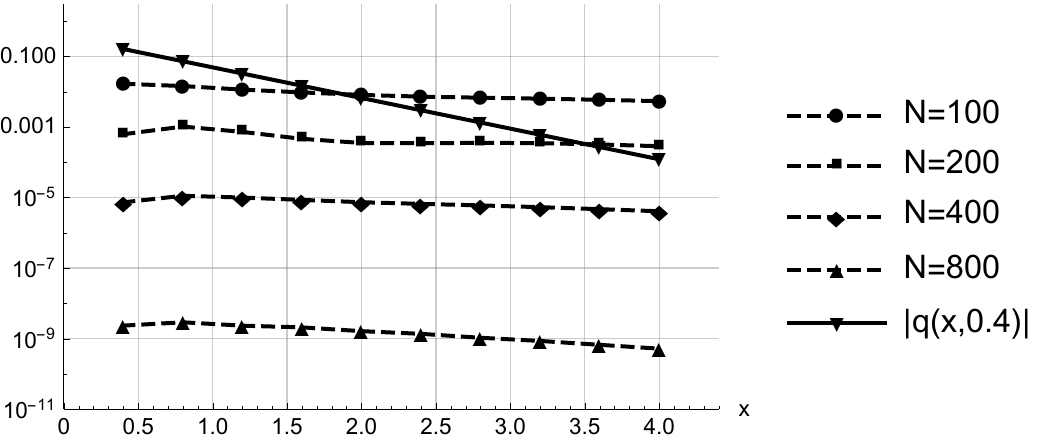}
  \caption{}
  \end{subfigure}
  \begin{subfigure}{0.45\textwidth}
  \includegraphics[width=\textwidth]{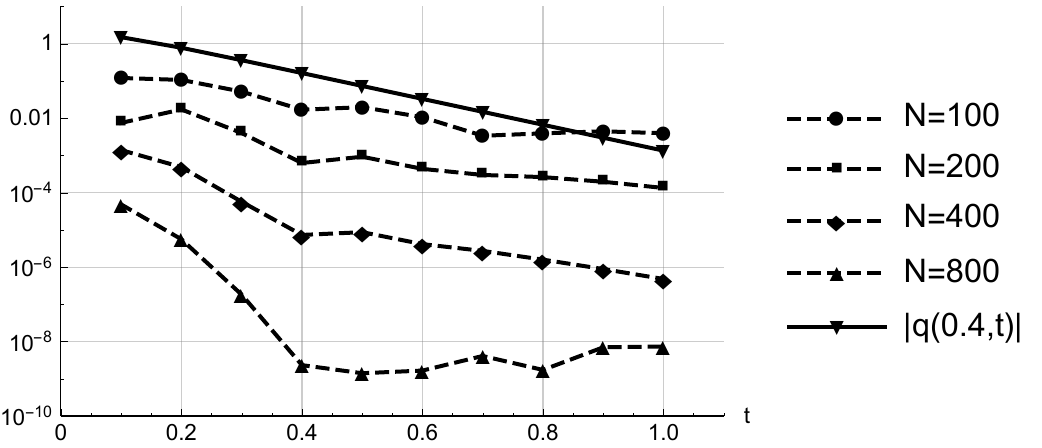}
  \caption{}
  \end{subfigure}
  \begin{subfigure}{0.45\textwidth}
  \includegraphics[width=\textwidth]{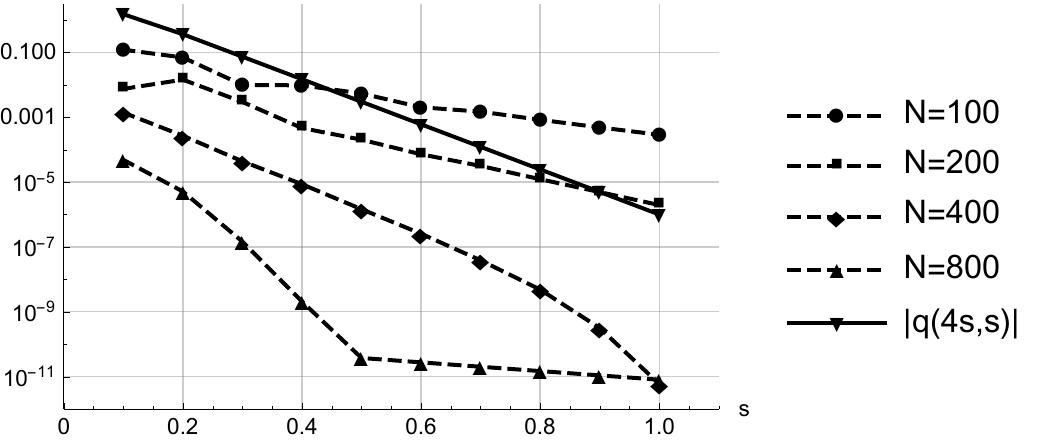}
  \caption{}
  \end{subfigure}

  \caption{The absolute value of the exact one-soliton solution (\ref{onesoliton}) with $\xi=1,\eta=1,x_0=0.4,\phi_0=0$ (solid lines) and the absolute errors of the numerical solution with different numbers of collocation points $N$ (dashed lines). Panel (a) shows the evaluations for $x\in [0.4,4],t=0.4$. Panel (b) shows the evaluations for $x=0.4,t\in [0.1,1]$. Panel (c) shows the evaluations for $x=4s,t=s, s \in [0.1,1]$. The tolerance of the error from computation of the spectral functions and contour truncations is set to $10^{-9}$.}
  \label{solitonerror}
\end{figure}

Figure \ref{solitonerror} shows the error plots of the solution using the NUTM with the initial and boundary values given by (\ref{onesoliton}) with $\xi=1,\eta=1,x_0=0.4$ along different lines in the $x,t$ quarter plane. The spectral convergence of the NUTM is demonstrated by the errors with fixed $x,t$ and varying $N$, the number of collocation points used in the solution of the RHP. The absolute errors are uniformly controlled and decreasing for increasing $x,t$ with fixed $N$. For increasing $t$, the NUTM even maintains relative accuracy with $N$ fixed, but sufficiently large. Although the absolute errors decrease exponentially for fixed $t$, the exponential decay rate of the solution is not captured exactly. As a result, the relative error with fixed $t$ increases as $x$ grows. This is expected, due to the fact that the jump matrix after deformation decays to the identity exponentially but not exactly at the decay rate of the solution. Indeed, how far we can deform the contour is restricted by the region of analyticity of $\gamma(k)$ and $\Gamma(k)$. For instance, when $q_0 \in \mathcal{S}_{2\alpha}$, $\gamma(k)$ is analytic within $-\alpha \leq \im{k}\leq \alpha$, along a horizontal segment of the contour $k=s+i\alpha$,
\[
e^{2i\theta(k;x,t)}=e^{-8 s \alpha t - 2 \alpha x + i (4 (s^2- \alpha^2) t + 2 s x)}.
\]
This is to be compared with the situation for the $1$-soliton solution (\ref{onesoliton}) which has an exponential decay rate $2\eta$ in the $x$ direction. The zeros of $a(k)$ are outside the strips $-\alpha \leq \im{k}\leq \alpha$, since $\abs{\alpha}<\abs{\eta\tanh(2\eta x_0)}<\abs{\eta}$. To capture the same exponential decay rate, a deformation of the horizontal contours up to the pole of  $\gamma(k)$ at $k=\xi+i \eta \tanh(2\eta x_0)$ is necessary. The restriction is not required if $t$ is sufficiently large, in which case the jump matrix along the deformed contour approaches the identity and is truncated before Im$(k)=\eta \tanh(2\eta x_0)$. For instance, consider the jump functions related to $P$ in Figure \ref{contour1}. When $t$ is sufficiently large, the jump function $\Delta P \Delta^{-1}$ in the top right corner of Figure \ref{contour1} is close to the identity matrix and is negligible. After the truncation, only the jump function $\Delta P J_1 \Delta^{-1}$ remains.

\subsection{Boundary conditions implicitly determined by given spectral data $A(k,T)$ and $B(k,T)$}
\subsubsection{A dressing argument}
If $A(k,T)$ and $B(k,T)$ are given directly, the deformation steps discussed in Section \ref{sec_linearizablebc} can still be performed provided that the deformations are within the regions where $A(k,T)$ and $B(k,T)$ are analytic. From the idea of the dressing method, as long as the spectral functions $A(k,T)$ and $B(k,T)$ satisfy (\ref{Sdetcondition}) and (\ref{globalrelation}), the RHP (\ref{rhp1}) generates solutions to the NLS equation.
\begin{proposition*}[The dressing method \cite{fokasbook}]
Suppose that the oriented smooth curve $\mathcal{L}$ divides the complex $k$-plane into the domains $D^+$ and $D^-$. Let $M^{+}(k;x,t)$ satisfy the following $2\times 2$ matrix RHP in the complex $k$-plane for all $x,t>0$,
\[
M^{+}(k;x,t)=M^{-}(k;x,t) e^{-i(kx+2k^2t)\hat{\sigma}_3} J(k), \, k\in \mathcal{L},
\]
where $J(k)$ is a $2\times 2$ unimodular matrix with $J_{11}=1$ or $J_{22}=1$ . Assume that the above RHP has a unique solution which is sufficiently smooth for all $x,t>0$. Define $Q(x,t)$ by
\[
Q(x,t)=i\lim_{k\rightarrow \infty} [\sigma_3, k M(x,t,k)],
\]
then $Q(x,t)$ satisfies the nonlinear equation
\[
iQ_t-Q_{xx}\sigma_3+2 Q^3 \sigma_3=0.
\]
\end{proposition*}
To ensure that (\ref{Sdetcondition}) is satisfied, we specify the ratio $h(k)=B(k)/A(k)$ for $k\in D_3$ since this is the quantity required in (\ref{rhp1}). The spectral data $A(k)$ and $B(k)$ are defined implicitly from (\ref{Sdetcondition}),
\[
1 + \lambda h(k) \overline{h(\overline{k})} =\frac{1}{A(k)\overline{A(\overline{k})}}, \,\,k\in \mathbb{R}\cup i\mathbb{R}.
\]
\begin{remark}
The global relation (\ref{globalrelationinf}) determines the value of $h(k)$ in the first quadrant. If $T=\infty$, $h(k)=b(k)/a(k)$ provided that $a(k)\neq 0$.
\end{remark}
\subsubsection{Numerical results}
\begin{figure}
  \centering
  \includegraphics[width=0.45\textwidth]{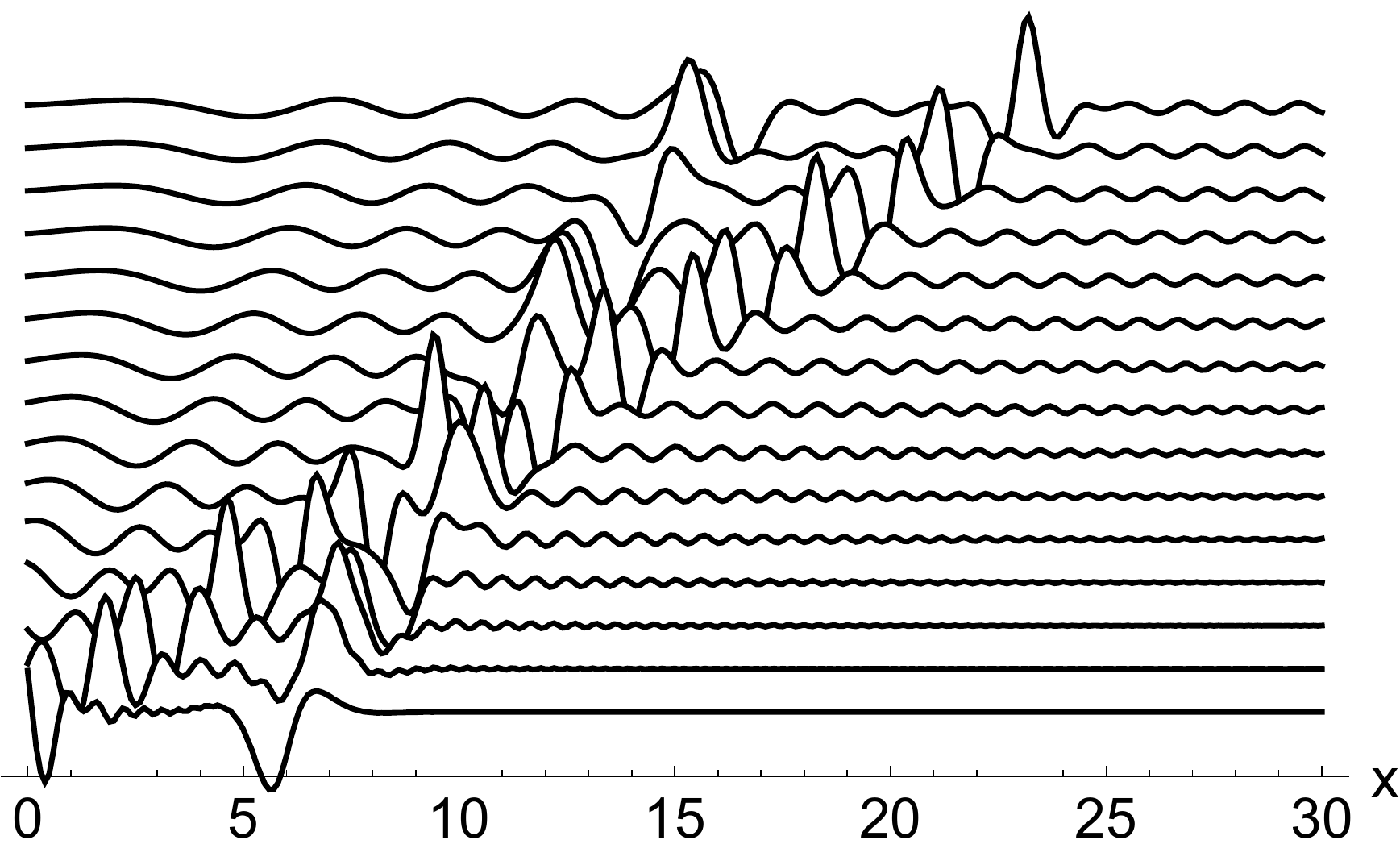}
  \includegraphics[width=0.45\textwidth]{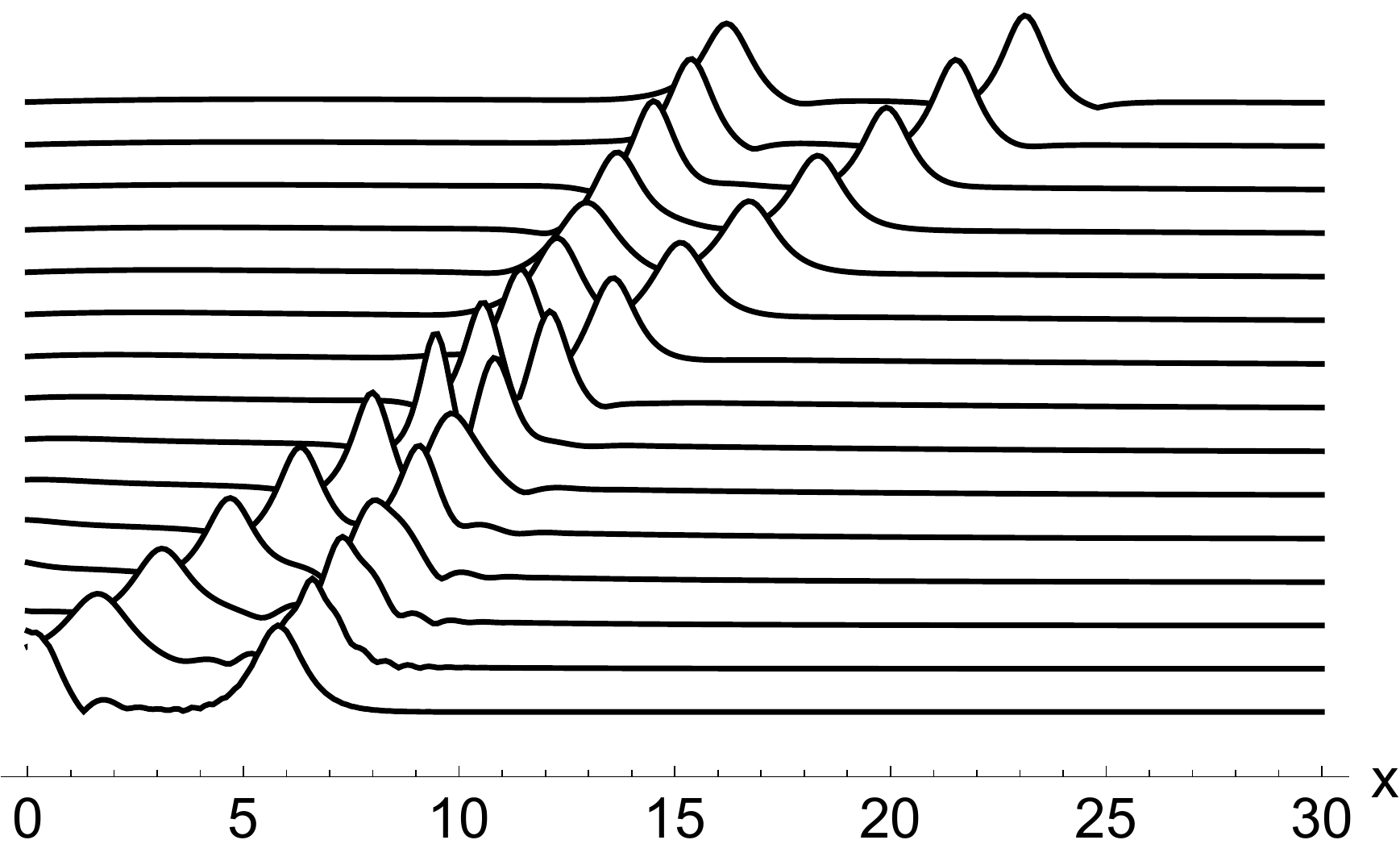}
  \caption{The numerical solution $q(x,t)$ on the domain $0\leq x\leq 30$ and $0.1\leq t\leq 3$ with spectral functions specified in Section 4.3. The solution contains two right-going solitons as well as dispersion. Left: The real part of $q(x,t)$. Right: The absolute value of $q(x,t)$. }\label{dressing1}
\end{figure}

\begin{figure}
  \centering
  \includegraphics[width=0.45\textwidth]{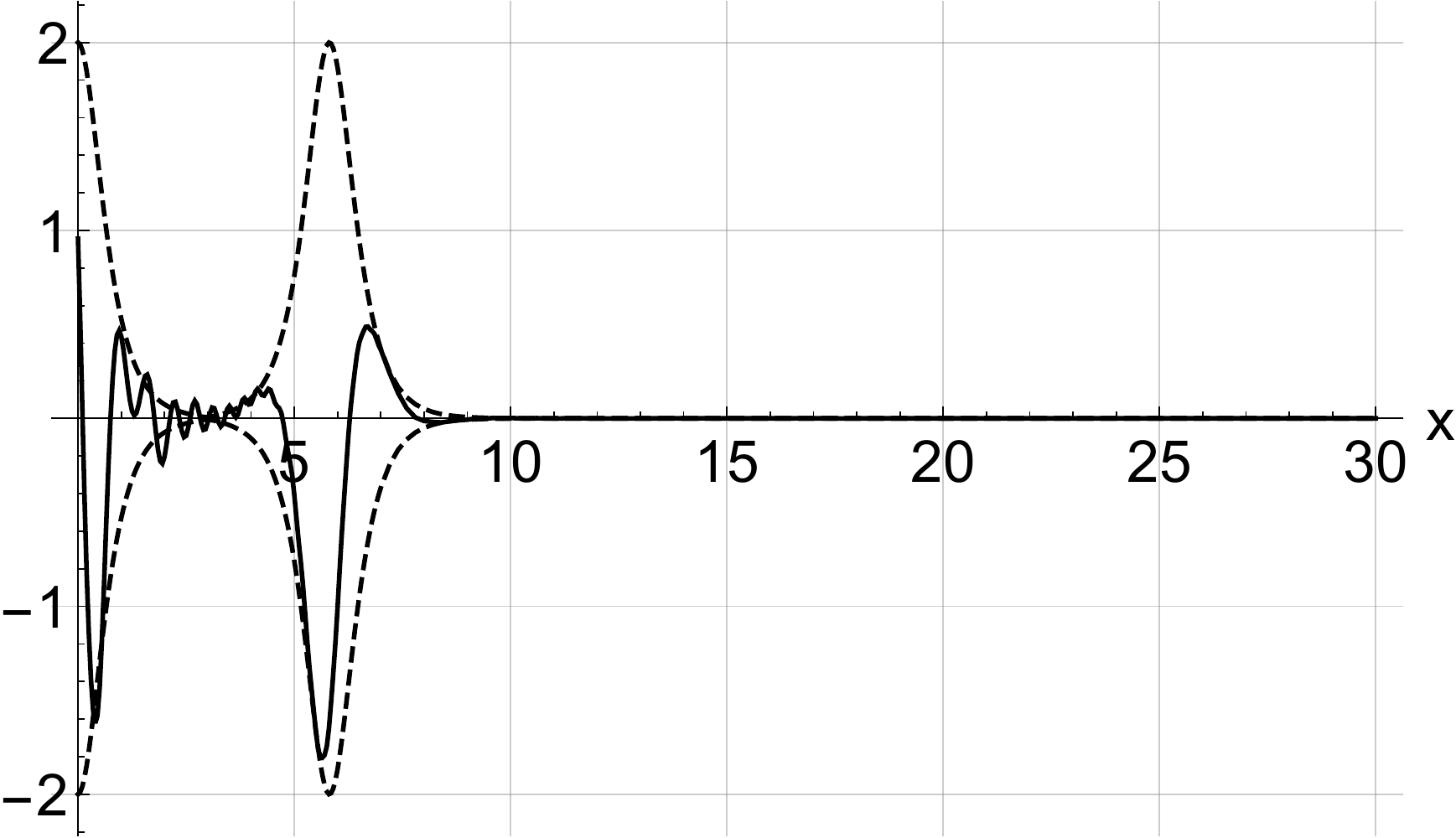}
  \includegraphics[width=0.45\textwidth]{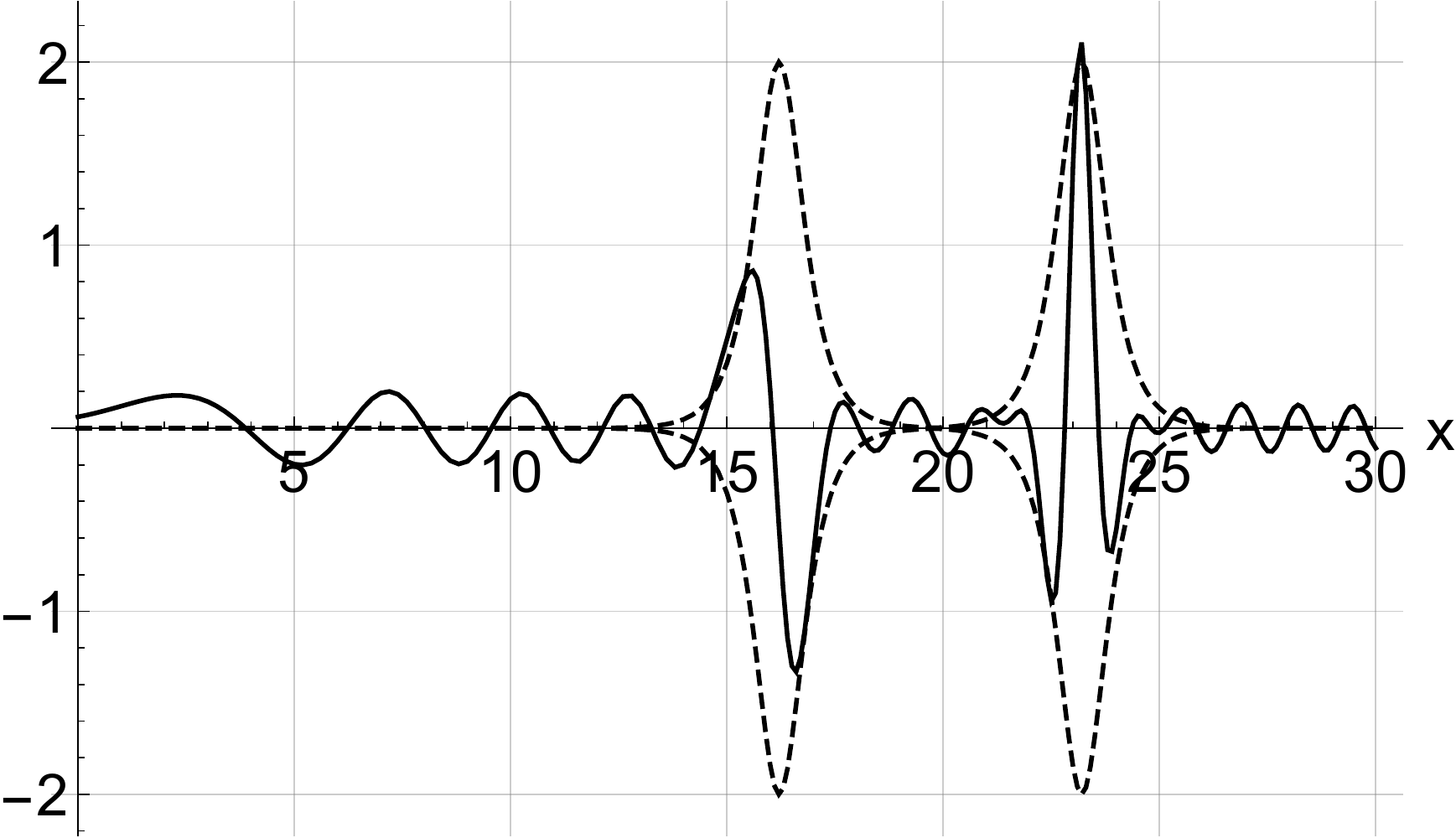}
  \caption{The real part of $q(x,t)$ (solid curve). The envelope of the exact right-going 2-soliton (dashed curve). Left: $q(x,0.1)$. Right: $q(x,2.9)$. }\label{dressing2}
\end{figure}

\begin{figure}
  \centering
  \begin{subfigure}{0.45\textwidth}
  \includegraphics[width=\textwidth]{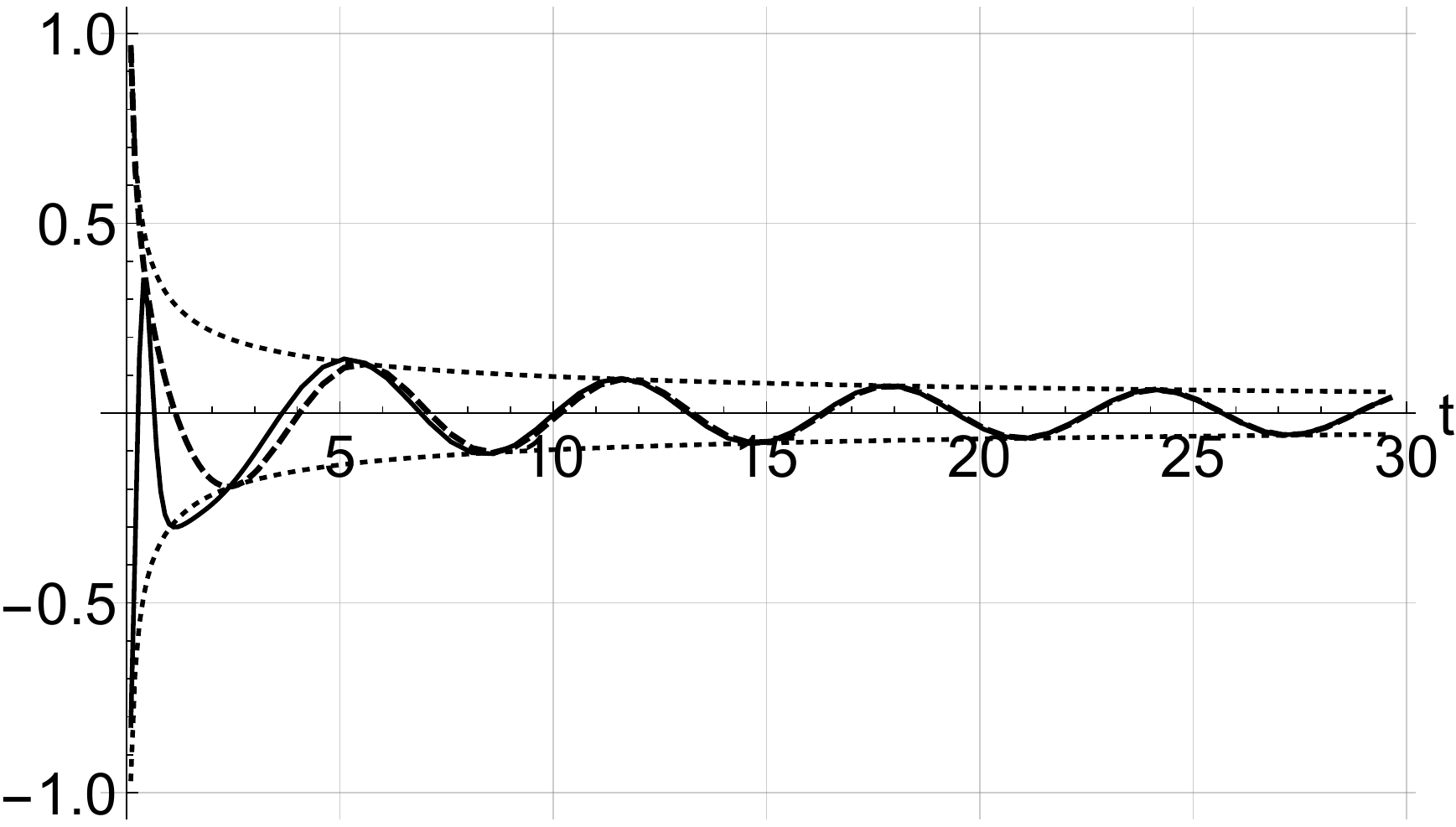}
  \caption{}
  \end{subfigure}
  \begin{subfigure}{0.45\textwidth}
  \includegraphics[width=\textwidth]{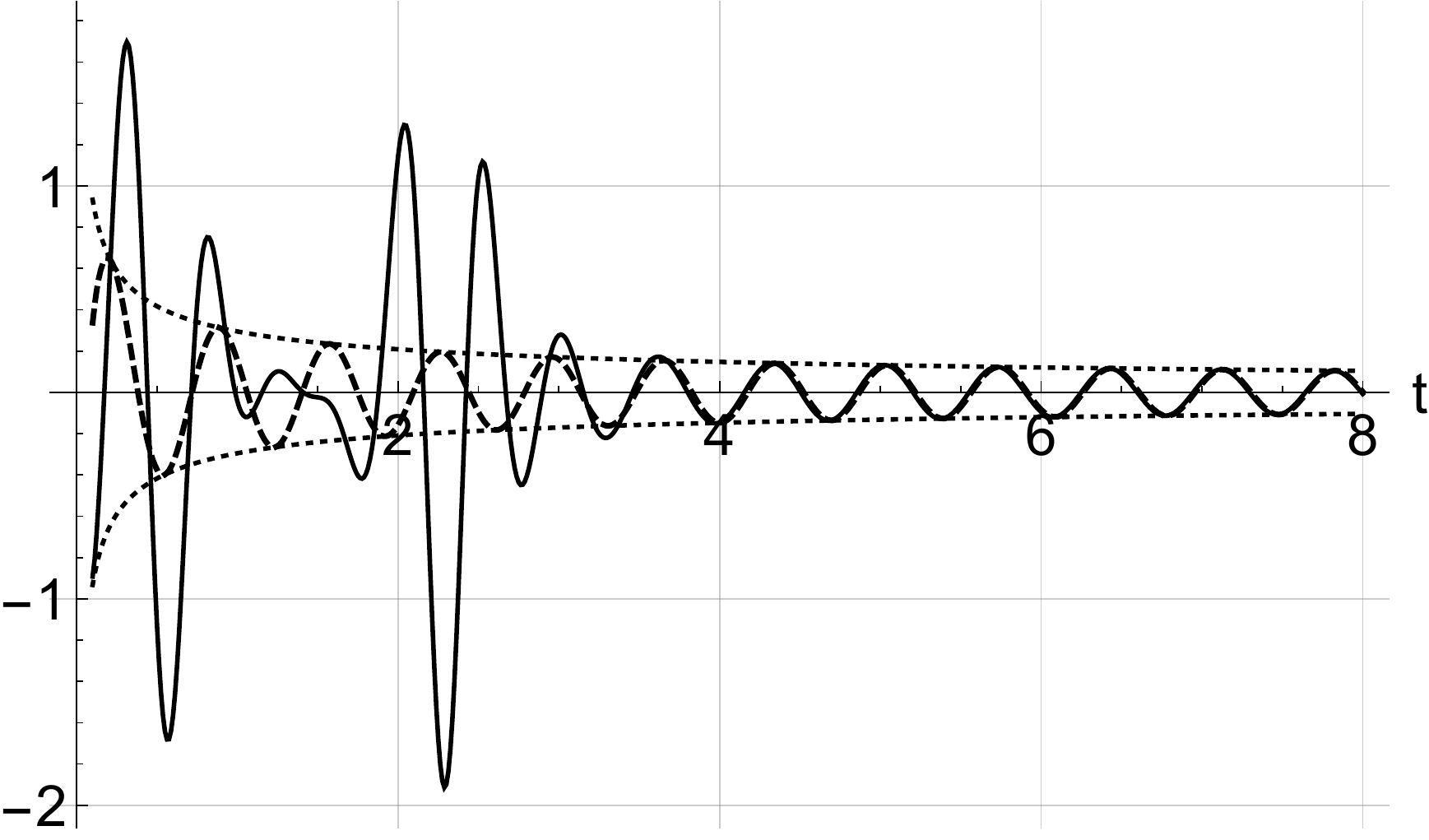}
  \caption{}
  \end{subfigure}\\
  \begin{subfigure}{0.45\textwidth}
  \includegraphics[width=\textwidth]{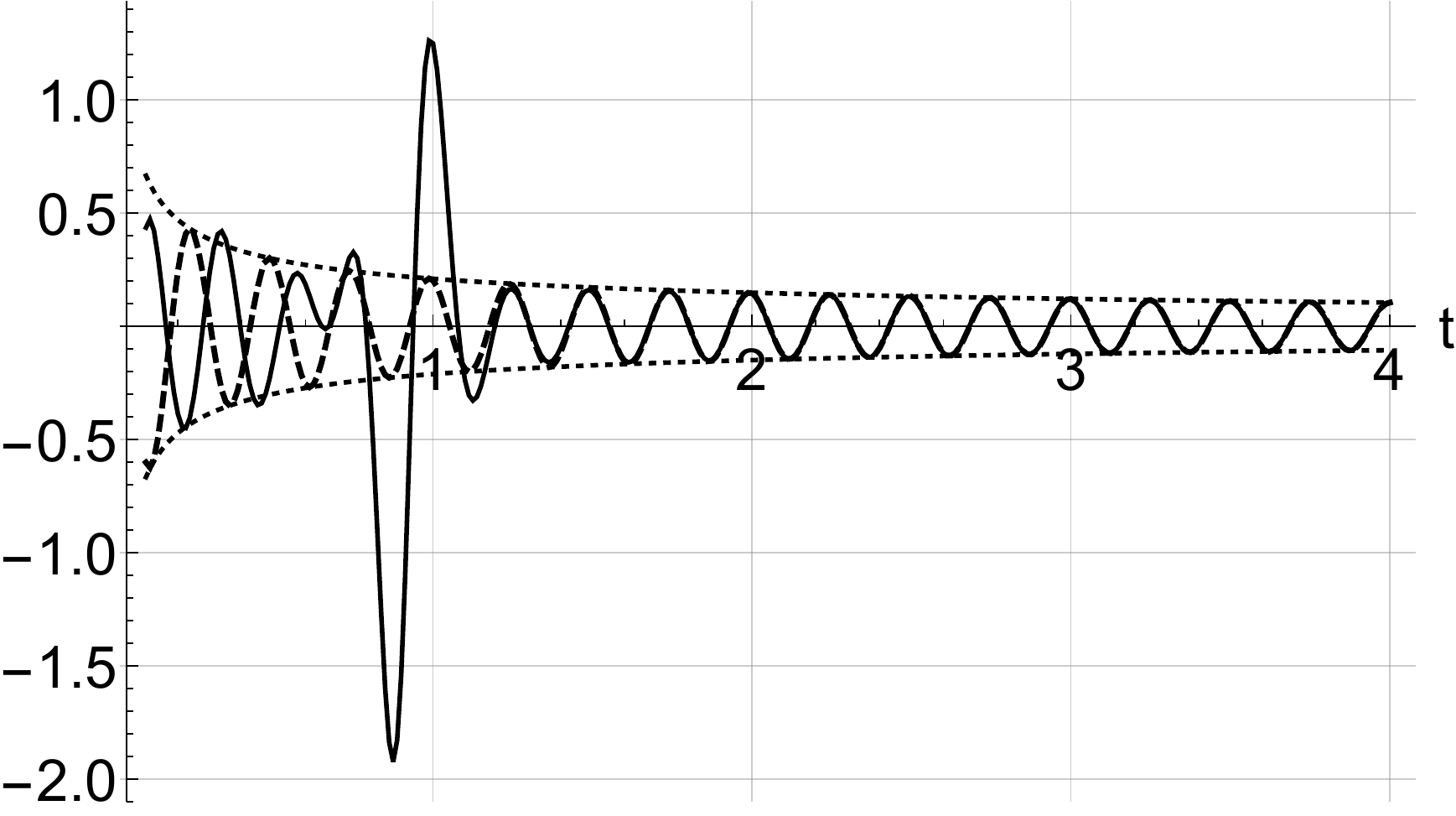}
  \caption{}
  \end{subfigure}
  \begin{subfigure}{0.45\textwidth}
  \includegraphics[width=\textwidth]{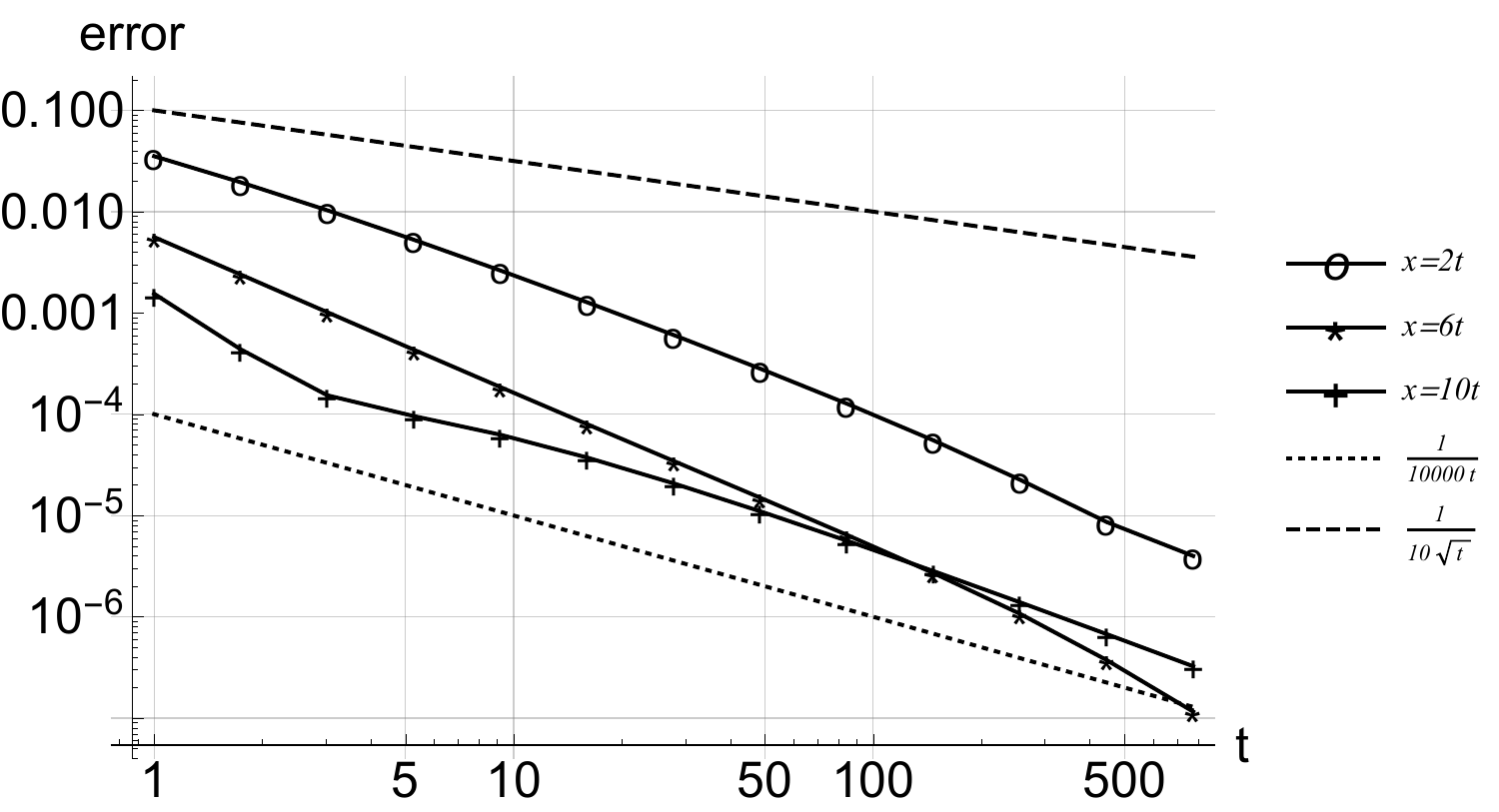}
  \caption{}
  \end{subfigure}
  \caption{The real part of the numerically computed $q(x,t)$ (solid curve), the real part of the dispersive wave from the asymptotic formula (\ref{asympformula}) (dashed curves) and its envelope (dotted curves). Panel (a) shows evaluation along $x=2t$. Panel (b) shows evaluation along $x=6t$. Panel (c) shows evaluation along $x=10t$. Panel (d) shows the errors compared with the asymptotic formula (\ref{asympformula}).}
  \label{dressing3}
\end{figure}

We solve the focusing NLS equation on the half-line with spectral functions $b(k)/a(k)=0$ and $B(k)/A(k)=1000k/\left(k-2(1+i)\right)^{5}$ for $\arg(k)\in [\pi,3\pi/2]$. In this case, $\gamma(k)=0$ and $\Gamma(k)=1000k/\left(k-2(1-i)\right)^{5}$. Furthermore, we impose two residue conditions at $k_1=-1+i$ and $k_2=-2+i$ with corresponding constants $c_1=100000,c_2=2$. The constant $1000$ in $\Gamma$ is chosen so that the dispersion is on the same order of the solitons for small $x,t$. The constants $c_1,c_2$ are chosen so that the interaction of the solitons is inside the domain. We plot the real part and the absolute value of the solution for $0\leq x\leq 20$ and $0.1\leq t\leq 3$ in Figure \ref{dressing1}. The solution contains two right-going solitons as well as dispersion. Two slices of the solution at $t=0.1$ and $t=2.9$ are shown in Figure \ref{dressing2}. We also observe that the soliton part of the numerical solution is similar to the exact two-soliton solution with its envelope plotted in dashed lines.
We compare our solution with the large $t$ asymptotics along $x/t=2,6,10$ in Figure \ref{dressing3}. Away from the solitons, the large $t$ asymptotics along $x/t=O(1)$ is described by (see \cite{fokas2005})
\begin{align}
q(x,t)=t^{-1/2}\alpha\left(-\frac{x}{4t}\right) \text{exp} \left(\frac{ix^2}{4t}+2i\lambda \alpha^2 \left(-\frac{x}{4t}\right)\log t + i \phi\left(-\frac{x}{4t}\right) \right) + o\left(t^{-\frac{1}{2}}\right),  \text{\quad as \quad} t\rightarrow \infty,
\label{asympformula}
\end{align}
where the amplitude $\alpha$ and the phase $\phi$ are given by
\begin{align*}
\alpha^2(k)=&\frac{\lambda}{4\pi} \log\left(1+\lambda \abs{\gamma(k)+\lambda \overline{\Gamma(k)}}^2\right),\\
\phi(k)=&6\lambda \alpha^2(k)\log 2 + \frac{\pi(2+\lambda)}{4}+ \arg\left(\gamma(k)+\lambda \overline{\Gamma(k)}\right)+ \arg \Gamma(-2i\lambda \alpha^2(k))+4\lambda \int_{-\infty}^{k} \log\abs{\mu-k} d\alpha^2(\mu).
\end{align*}
In practice, to avoid computing the integral in the formula for the phase $\phi$, since $\phi$ is constant with fixed $x/t$, we choose it so that the errors in Figure (\ref{dressing3}d) show a trend of decreasing errors with order $O(1/t)$.

\section{Using large $k$ expansions for computing $q(x,t)$ for small $x,t$}
\label{sec_fasterdecay}
For a problem on the whole line, if the initial condition is in Schwartz class, the reflection coefficient is also in Schwartz class \cite{beals1984}. Therefore no deformation for the associated RHP is required when $x,t$ are small because the jump matrix decays rapidly to the identity matrix when $\abs{k}$ becomes large. However, this is not true for the RHP (\ref{rhp1}) from problems on the half-line, or even for integrals that arise in the linear case.
\begin{example}
Consider the NLS equation with the homogeneous Neumann boundary condition $g_1(t)=0$ and initial condition $q_0(x)=e^{-x^2}+ i \text{sech}^2(x)$. This is a linearizable boundary condition and we can solve for the spectral functions $B(k,\infty)$, $A(k,\infty)$ using symmetries of the global relation $k\rightarrow -k$. The associated RHP (\ref{rhp5})
 for $\Phi(x,t,k)$ is formulated with the jump condition on the cross $k\in \mathbb{R}\cup i\mathbb{R}$ in Figure \ref{rhp1contour},
\begin{align}
\Phi_{+}(k;x,t)=\Phi_{-}(k;x,t)J(k;x,t).
\label{rhp5}
\end{align}
The jump matrices are the same as in (\ref{rhp1jump}) except that (\ref{rhp1uppergamma}) is replaced by
\begin{align*}
\Gamma(k)=\frac{\lambda \overline{b(-\overline{k})}}{a(k)\Delta_1(k)} &, ~  \text{arg }k\in [0,\pi],
\end{align*}
where
\begin{align*}
\Delta_1(k)=a(k)  \overline{a(-\overline{k})} -\lambda b(k)
\overline{b(-\overline{k})}&, ~\text{arg } k \in [0,\pi].
\end{align*}
The functions $\gamma(k)$ and $\Gamma(k)$ are $O(1/k)$ as $k\rightarrow \infty$, see Figure \ref{gammaplot}.
\begin{figure}
  \makebox[\textwidth][c]{
 \begin{subfigure}{0.3\textwidth}
  \includegraphics[width=1.\textwidth]{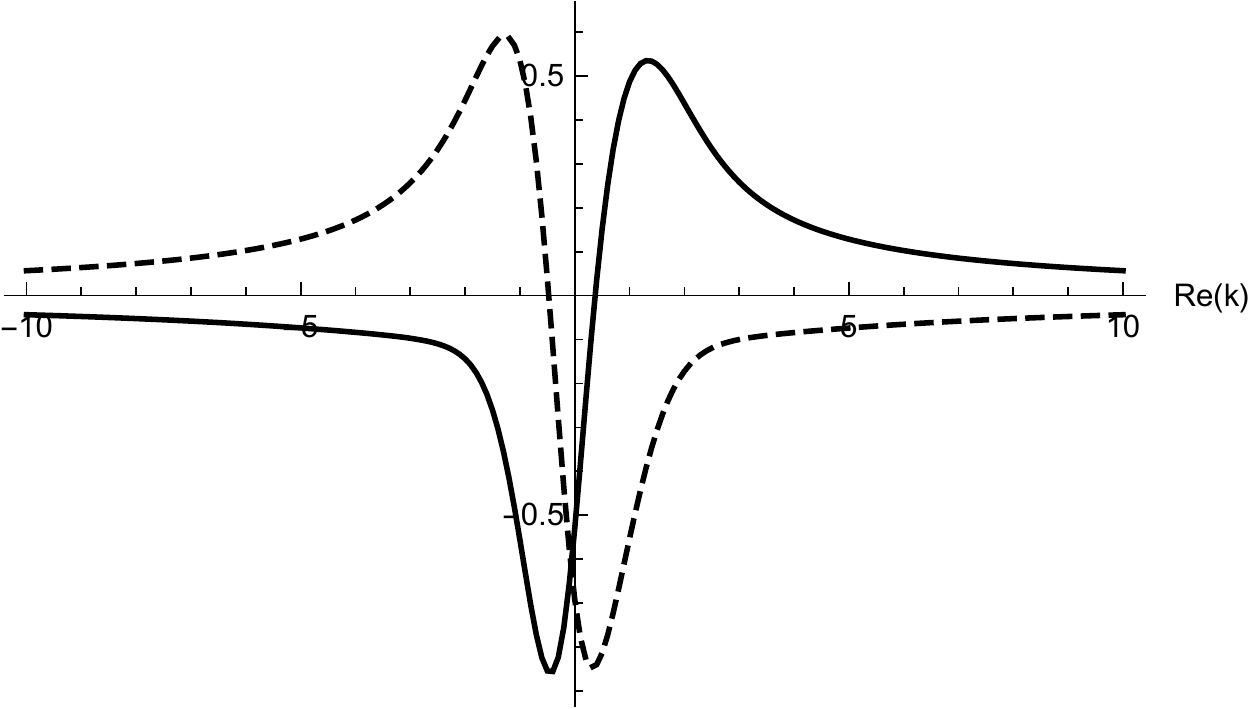}
\end{subfigure}
 \begin{subfigure}{0.3\textwidth}
    \includegraphics[width=1.\textwidth]{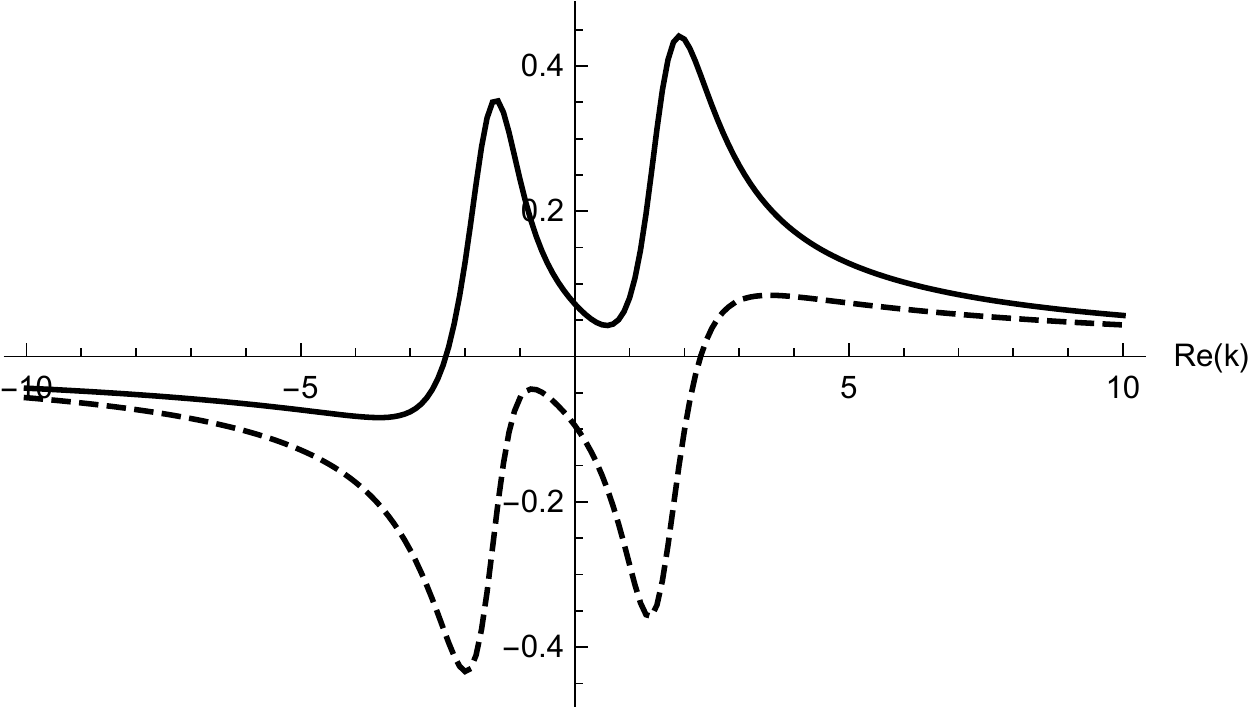}
\end{subfigure}
 \begin{subfigure}{0.3\textwidth}
    \includegraphics[width=1.\textwidth]{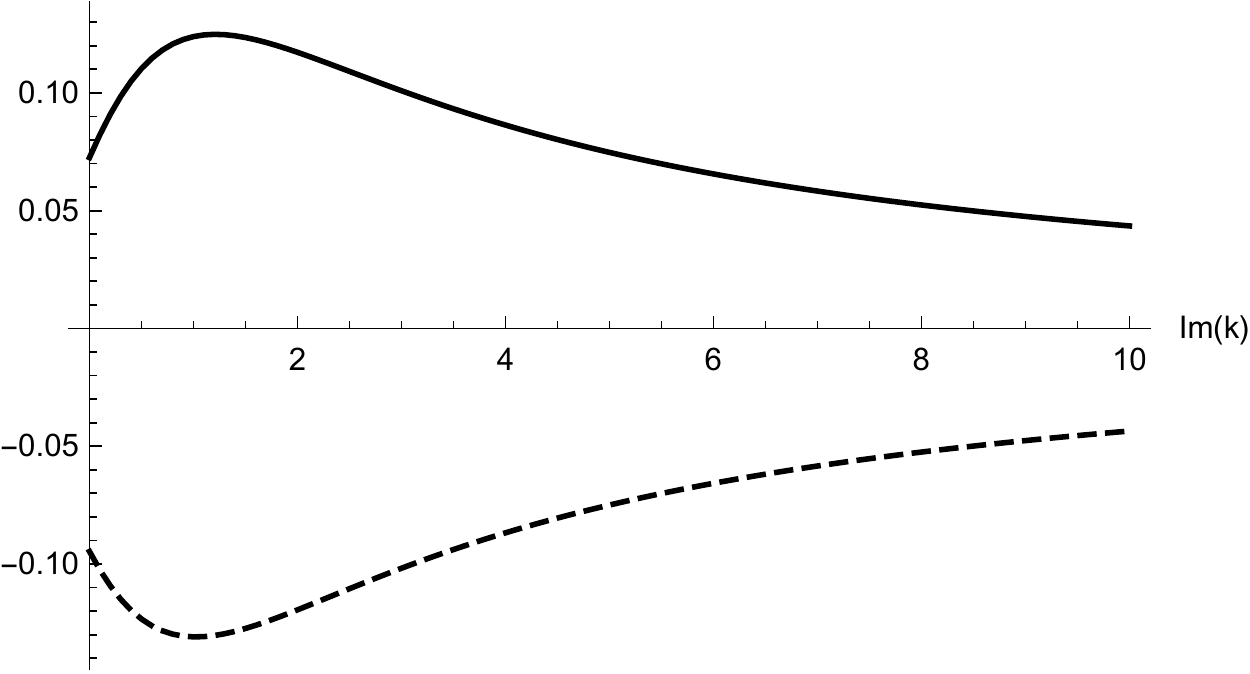}
\end{subfigure}
  }
  \caption{The plots of $\gamma(k),\Gamma(k)$ along the jump contour. Left: $\gamma(k)$ on the real axis. Middle: $\Gamma(k)$ on the real axis. Right: $\Gamma(k)$ on the positive imaginary axis. The real parts are plotted with solid curves and the imaginary parts are plotted with dashed curves. The initial condition is $q_0(x)=e^{-x^2}+i\text{ sech}^2(x)$ and the boundary condition is $g_1(t)=0$. The functions $\gamma$ and $\Gamma$ are of order $1/k$.}
  \label{gammaplot}
\end{figure}
\end{example}

\begin{figure}
  \centering
  \includegraphics[width=0.6\textwidth]{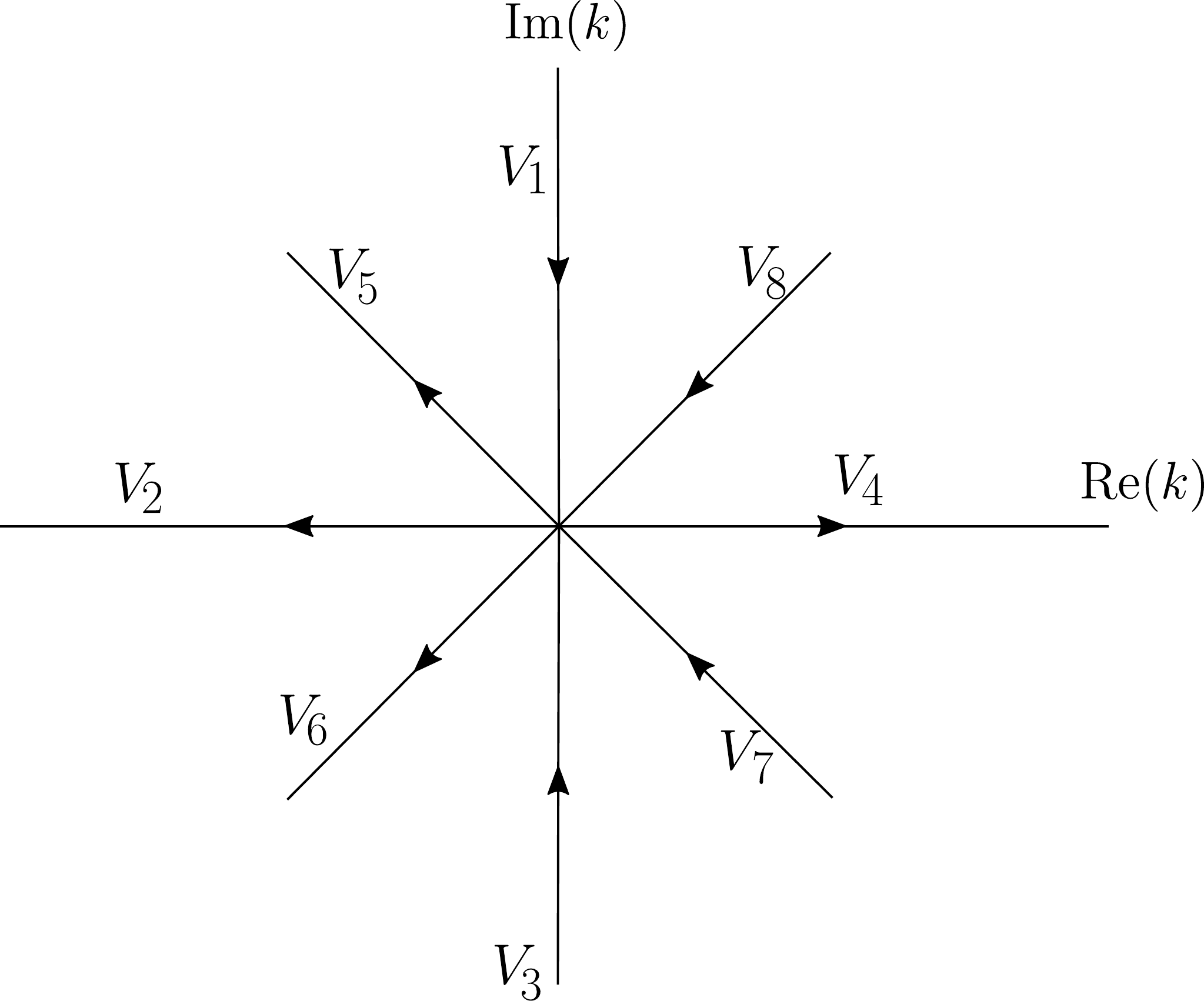}\\
  \caption{The contour of RHP (\ref{rhp5}). Jumps on $V_5,V_6,V_7,V_8$ are introduced to improve the rate at which the jump matrices on $V_1,V_2,V_3,V_4$ approach the identity matrix. }\label{rhp5contour}
\end{figure}

\begin{figure}
  \centering
  \includegraphics[width=0.6\textwidth]{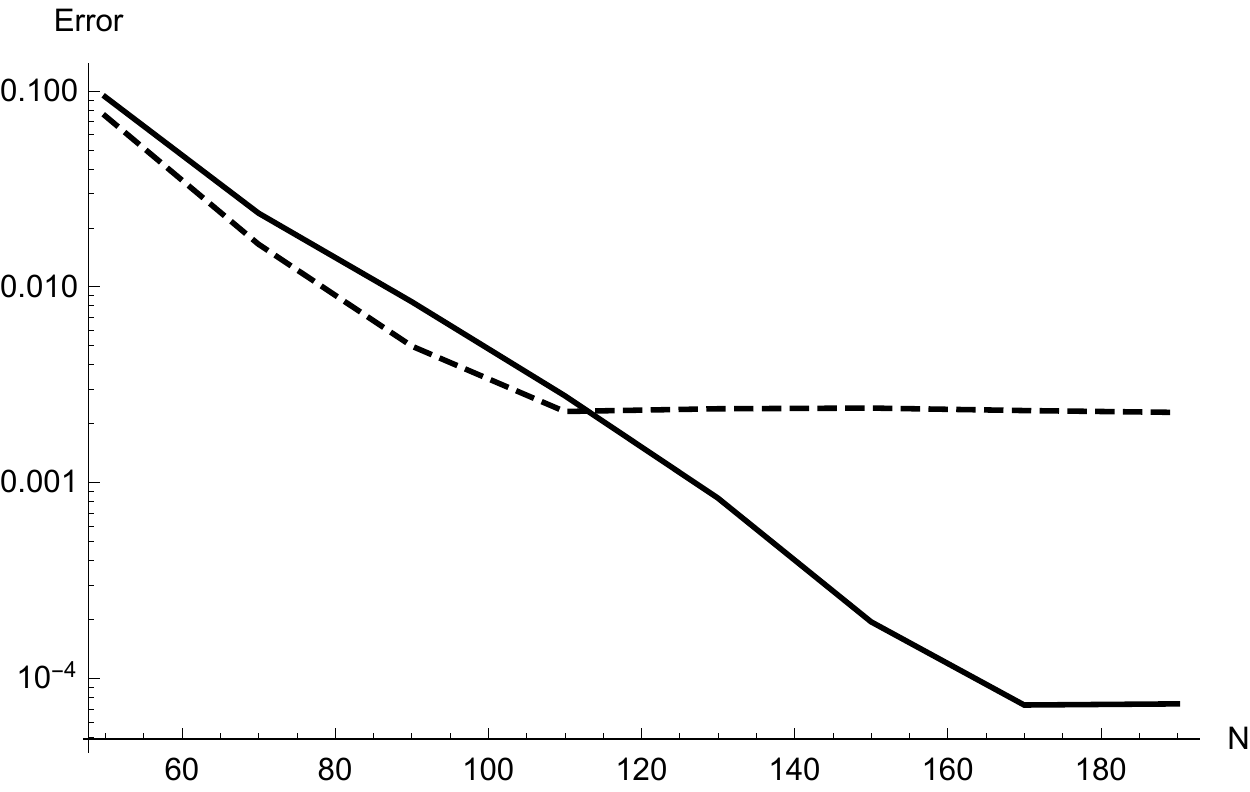}\\
  \caption{The absolute error for computing $q(0.5,0)$ with different number of collocation points $N$. The dashed line is computed using the undeformed contour in Figure \ref{rhp1contour}. The solid line is computed using the contour in Figure \ref{rhp5contour}. The flattening in both curves is due to the truncation error.}\label{accuracyplot1}
\end{figure}

With the large $k$ expansions derived in the Appendix, we can set up RHPs with jump functions that tend to the identity matrix faster. We define
\begin{align*}
\gamma_0(k) =  \frac{q_0(0)}{2i\left(k-\hat{k}\right)}, \quad
\Gamma_0(k) =  \frac{\lambda \overline{q_0(0)}}{2i\left(k-\hat{k}\right)},
\end{align*}
\begin{align*}
\gamma(k) = \gamma_0(k) + \gamma_r(k), \quad
\Gamma(k) = \Gamma_0(k) + \Gamma_r(k),
\end{align*}
so that $\gamma_0(k)$ and $\Gamma_0(k)$ have the same large $k$ behavior as $\gamma(k)$ and $\Gamma(k)$ to the leading order. To avoid introducing unnecessary residue conditions, we choose $\hat{k}=1-2i$. After separating the $O(1/k)$ terms, we get a RHP which has the jump matrix approaches the identity matrix with $O(1/k^2)$ as $k\rightarrow \infty$ shown in Figure \ref{rhp5contour}. The RHP for $\Phi(x,t,k)$ is formulated with a jump condition on the eight rays through the origin $\{r=\rho e^{i s } |  \rho\in [0,+\infty),\,\, s=0,\pi/4,2\pi/4,3\pi/4,4\pi/4,5\pi/4,6\pi/4,7\pi/4\}$.
and $\Phi(x,t,k)$ satisfies
\begin{align}
\Phi_{+}(k;x,t)=\Phi_{-}(k;x,t)V(k;x,t),
\label{rhp5}
\end{align}
with jump matrices
\begin{align}
V(k;x,t)= \begin{cases}
                        V_4(k;x,t), & \mbox{arg } k = 0, \\
                        V_1(k;x,t), & \mbox{arg } k = \frac{\pi}{2},\\
                        V_2(k;x,t), & \mbox{arg } k = \pi,\\
                        V_3(k;x,t), & \mbox{arg } k = \frac{3\pi}{2},\\
                        V_5(k;x,t), & \mbox{arg } k = \frac{3\pi}{4}, \\
                        V_6(k;x,t), & \mbox{arg } k = \frac{5\pi}{4},\\
                        V_7(k;x,t), & \mbox{arg } k = \frac{7\pi}{4},\\
                        V_8(k;x,t), & \mbox{arg } k = \frac{\pi}{4},
    \end{cases}
\end{align}
where
\begin{align*}
V_1(k;x,t)= & \left[\begin{matrix}
1 & 0 \\
-\Gamma_r(k)e^{2i\theta(k;x,t)} & 1
\end{matrix}\right], \quad
V_3(k;x,t)=\left[\begin{matrix}
1 & -\lambda \overline{\Gamma_r(\overline{k})} e^{-2i\theta(k;x,t)} \\
0 & 1
\end{matrix}\right], \\
V_4(k;x,t)=&\left[\begin{matrix}
1+\lambda  \overline{\gamma_r(\overline{k})} \gamma_r(k) & \gamma_r(k) e^{-2i\theta(k;x,t)} \\
\lambda \overline{\gamma_r(\overline{k})} e^{2i\theta(k;x,t)} & 1
\end{matrix}\right],
\\
V_2(k;x,t) =& \left[\begin{matrix}
1  & -\left(\lambda \overline{\Gamma_r(\overline{k})}  + \gamma_r(k) \right)e^{-2i\theta(k;x,t)} \\
-\left( \lambda \overline{\gamma_r(\overline{k})} +  \Gamma_r(k) \right)e^{2i\theta(k;x,t)} &
 1 +\left( \lambda \overline{\Gamma_r(\overline{k})} +  \gamma_r(k) \right)\left( \lambda \overline{\gamma_r(\overline{k})} + \Gamma_r(k)\right)
\end{matrix}\right],\\
V_5(k;x,t)= & \left[\begin{matrix}
1 & -\left( \lambda \overline{\gamma_0(\overline{k})} +  \Gamma_0(k) \right) e^{-2i\theta(k;x,t)} \\
0 & 1
\end{matrix}\right], \quad
V_6(k;x,t)=  \left[\begin{matrix}
1 & 0 \\
- \left(\lambda \overline{\Gamma_0(\overline{k})}  + \gamma_0(k) \right) e^{2i\theta(k;x,t)} & 1
\end{matrix}\right], \\
V_7(k;x,t)=&\left[\begin{matrix}
1 & -\left(\lambda \overline{\Gamma_0(\overline{k})} + \gamma_0(k) \right)e^{-2i\theta(k;x,t)} \\
0 & 1
\end{matrix}\right],
 \quad
V_8(k;x,t)= \left[\begin{matrix}
1 & 0 \\
-\left(\Gamma_0(k)+\lambda \overline{\gamma_0(\overline{k})} \right)e^{2i\theta(k;x,t)} & 1
\end{matrix}\right].
\end{align*}
Figure \ref{accuracyplot1} shows the log-linear plot of the absolute error for computing $q(0.5,0)$ with different numbers of collocation points $N$. In the computation, the segments of contours are truncated when the jump matrix is close to the identity matrix $\norm{V_m(k)-I}_2<10^{-8}$, $m=1,2,\dots,8,$ or when the contour reaches a large circle centered at the origin with radius $50$. In Figure \ref{accuracyplot1}, the dashed curve shows the error computed using the undeformed contour and the solid curve shows the error computed using the contour in Figure~\ref{rhp5contour}. Both curves decay exponentially when $N$ is not large. Since the jump matrix of RHP (\ref{rhp5}) approaches the identity matrix faster, the flattening in the solid curve appears later than the dashed curve. It is possible to perform the asymptotic analysis for higher order terms in the previous section and remove more terms so that the decay of the jump matrix is faster than $O(1/k^2)$. The calculation starts to become lengthy, however.
\begin{remark}
For linearizable boundary conditions, if the solution of the half-line problem can be mapped to a smooth solution of the whole-line problem, the jump matrix $J_2$ in (\ref{rhp1}) actually decays exponentially to the identity matrix. Therefore there is no need to introduce modifications to $V_2$, $V_5$ and $V_6$ in (\ref{rhp5}). In this case, $\Gamma(k)$ is automatically analytically extended to the first quadrant by the analyticity of $b(k)$ and $a(k)$. Then it is possible to deform the jump contour $J_1$ and $J_3$ in (\ref{rhp1}) to the positive real line on top of $J_4$ and this new RHP is the same as the RHP in the whole-line problem with the initial values on the negative real line defined properly corresponding to the boundary condition using symmetry.
\end{remark}

\section{Acknowledgements}
The authors gratefully acknowledge support from the US National Science Foundation under grants NSF-DMS-1522677 (BD,XY) and NSF-DMS-1945652 (TT). Any opinions, findings, and conclusions or recommendations expressed in this material are those of the authors and do not necessarily reflect the views of the funding sources.

\section*{Appendix: Large $k$ asymptotics of the spectral functions}
We consider the large $k$ asymptotics of the spectral functions $a(k),b(k)$, $A(k,\infty),B(k,\infty)$ and $A(k,T),B(k,T)$. This is used in Section \ref{sec_fasterdecay} to improve the decay of the jump matrix. It is known that the global relation characterizes the initial and boundary data of a solution to the NLS equation while the large $k$ asymptotics of the global relation characterizes the compatibility of initial and boundary data at $x=0$ and $t=0$ \cite{huang2018}. Some of these asymptotic results are found in \cite{fokas2012}. We re-examine these results with more terms, comparing them with examples from the LS equation, and use them to improve the decay of the jump matrix for the RHP (\ref{rhp1}). The initial data $q_0(x)$ and the boundary data $g_0(t)$, $g_1(t)$ are assumed to have sufficient smoothness and decay at infinity so that all relevant integrals are well defined.
\subsection*{Large $k$ asymptotics of $a(k),b(k)$}
Recall that the spectral functions $a(k)=\phi_2(0,k)$ and $b(k)=\phi_1(0,k)$ are defined using the linear Volterra integral equations (\ref{phi1}) and (\ref{phi2}).
We use the following expansions for large $k$ in the upper half-plane $\im{k}\geq 0$,
\alpheqn
\begin{align}
\label{phi1a} \phi_1(x,k)&=\frac{f_{11}(x)}{k}+\frac{f_{12}(x)}{k^2}+O\left(\frac{1}{k^3}\right),\\
\label{phi2a} \phi_2(x,k)&=1+\frac{f_{21}(x)}{k}+O\left(\frac{1}{k^2}\right).
\end{align}
\resetalpheqn
Substituting the expansions into (\ref{phi1}) and (\ref{phi2}) and matching terms with different powers of $k$, we obtain the following:
\begin{itemize}
\item Using
\begin{align*}
  -\int_x^{\infty} e^{-2ik(x-y)}q_0(y) \left(1+\frac{f_{21}(x)}{k}+O\left(\frac{1}{k^2}\right)\right) dy &= -\int_x^{\infty} e^{-2ik(x-y)}q_0(y)dy + O\left(\frac{1}{k^2}\right)\\
 &= \frac{q_0(x)}{2ik}+\frac{1}{2ik} \int_x^{\infty} e^{-2ik(x-y)}q'_0(y)dy +O\left(\frac{1}{k^2}\right),
\end{align*}
we obtain the coefficient of the $O(1/k)$ term in (\ref{phi1a}),
\[
f_{11}(x)=\frac{q_0(x)}{2i}.
\]
\item From
\begin{align*}
 &-\lambda \int_x^{\infty} \bar{q}_0(y) \phi_1(y,k)dy + O\left(\frac{1}{k^2}\right) = -\frac{\lambda}{2i k} \int_x^{\infty} \abs{q_0(y)}^2 dy +O\left(\frac{1}{k^2}\right),
\end{align*}
we get the coefficient of the $O(1/k)$ term in (\ref{phi2a}),
\[f_{21}(x)= -\frac{\lambda}{2i} \int_x^{\infty} \abs{q_0(y)}^2 dy.\]
\item Using
\begin{align*}
 &\frac{1}{2ik} \int_x^{\infty} e^{-2ik(x-y)}q'_0(y)dy -\int_x^{\infty} e^{-2ik(x-y)}q_0(y) \left( -\frac{\lambda}{2i k} \int_y^{\infty} \abs{q_0(z)}^2 dz \right) dy +O\left(\frac{1}{k^3}\right)\\
=& -\frac{q'_0(x)}{(2ik)^2} - \frac{ \lambda q_0(x)   \int_x^{\infty} \abs{q_0(z)}^2 dz }{(2ik)^2}\\
&\quad +\frac{\lambda}{(2ik)^2} \int_x^{\infty} e^{-2ik(x-y)} \left(  q_0(y) \abs{q_0(y)}^2 - q'_0(y) \int_y^{\infty} \abs{q_0(z)}^2 dz   \right)dy +O\left(\frac{1}{k^3}\right)\\
=& -\frac{q'_0(x)}{(2ik)^2} - \frac{ \lambda q_0(x)   \int_x^{\infty} \abs{q_0(y)}^2 dy }{(2ik)^2}+O\left(\frac{1}{k^3}\right),
\end{align*}
we get the coefficient of the $O(1/k^2)$ term in (\ref{phi1a})
\[
f_{12}(x)=-\frac{q'_0(x)}{(2i)^2} - \frac{ \lambda q_0(x)   \int_x^{\infty} \abs{q_0(y)}^2 dy }{(2i)^2}.
\]
\end{itemize}
As a result, we obtain expansions for $a(k),b(k)$:
\begin{align*}
a(k)&=\phi_2(0,k)= 1 -\frac{\lambda}{2i k} \int_0^{\infty} \abs{q_0(y)}^2 dy +O\left(\frac{1}{k^2}\right), \\
b(k)&=\phi_1(0,k)= \frac{q_0(0)}{2ik} -\frac{q'_0(0) + \lambda q_0(0)   \int_0^{\infty} \abs{q_0(y)}^2 dy }{(2ik)^2} +O\left(\frac{1}{k^3}\right).
\end{align*}
\subsection*{Large $k$ asymptotics of $A(k,T),B(k,T)$}
For $T<\infty$, the spectral functions $\overline{A(\overline{k},T)}=\Phi_2(T,k)$ and $-e^{-4ik^2T}B(k,T)=\Phi_1(T,k)$ are defined using the linear Volterra integral equations (\ref{phi3}) and (\ref{phi4}).
We use the following expansions for large $k\in \mathbb{C}$,
\alpheqn
\begin{align}
\label{phi3a} \Phi_1(t,k)&=\frac{h_{11}(t)}{k}+\frac{\hat{h}_{11}(t)e^{-4ik^2 t}}{k} +\frac{h_{12}(t)}{k^2}+\frac{\hat{h}_{12}(t)e^{-4ik^2 t}}{k^2}+O\left(\frac{1}{k^3}\right)+O\left(\frac{e^{-4ik^2t }}{k^3}\right),\\
\label{phi4a}
\Phi_2(t,k)&=1+\frac{h_{21}(t)}{k}+\frac{\hat{h}_{22}(t)e^{-4ik^2t}}{k^2}+\frac{\hat{h}_{23}(t)e^{-4ik^2t}}{k^3}+O\left(\frac{1}{k^2}\right)+O\left(\frac{e^{-4ik^2 t}}{k^4}\right),
\end{align}
\resetalpheqn
where the terms depending on $e^{-4ik^2t}$ are separated. Substituting the expansions into (\ref{phi3}) and (\ref{phi4}) and matching terms with different powers of $k$, we obtain the following:
\begin{itemize}
\item We have
\begin{align*}
&\int_0^{t} 2k\lambda \overline{g_0(\tau)} \left(\frac{h_{11}(\tau)}{k}+ \frac{\hat{h}_{11}(\tau)e^{-4ik^2 \tau}}{k} \right)+ i\lambda \abs{g_0(\tau)^2}d\tau + O\left(\frac{1}{k}\right)+O\left(\frac{e^{-4ik^2t}}{k^3}\right)\\
= & \lambda \int_0^{t} 2\overline{g_0(\tau)} h_{11}(\tau) + i \abs{g_0(\tau)^2}d\tau + \int_0^{t} 2\lambda \overline{g_0(\tau)} \hat{h}_{11}(\tau)e^{-4ik^2 \tau}d\tau \\
&+O\left(\frac{1}{k}\right)+O\left(\frac{e^{-4ik^2t}}{k^3}\right)
\\
=&  \lambda \int_0^{t} 2\overline{g_0(\tau)} h_{11}(\tau) + i \abs{g_0(\tau)^2}d\tau \\
&+ \frac{2i\lambda}{4k^2}\left(\overline{g_0(t)} \hat{h}_{11}(t)e^{-4ik^2t} - \overline{g_0(0)} \hat{h}_{11}(0) -   \int_0^{t}  e^{-4ik^2 \tau} \frac{d}{d\tau} \left( \overline{g_0(\tau)} \hat{h}_{11}(\tau) \right) d\tau \right)\\
&+ O\left(\frac{1}{k}\right)+O\left(\frac{e^{-4ik^2t}}{k^3}\right).
\end{align*}
Comparing this with $O(1)+O(e^{-4ik^2 t}/k^2)$ terms in (\ref{phi4}), we obtain
\begin{align*}
2\overline{g_0(\tau)} h_{11}(\tau) + i \abs{g_0(\tau)^2}=0 \quad \Rightarrow \quad h_{11}(t) = \frac{g_0(t)}{2i},\\
\end{align*}
and
\begin{align*}
 \hat{h}_{22}(t)=\frac{i\lambda}{2}\left(\overline{g_0(t)} \hat{h}_{11}(t)\right).
\end{align*}
\item Consider
\begin{align*}
& \int_0^{t} e^{-4ik^2(t-\tau)} \left( -i \lambda \abs{g_0^2} \left(\frac{h_{11}}{k}+\frac{\hat{h}_{11}e^{-4ik^2 \tau}}{k}\right)  \right)  d\tau \\
&+\int_0^{t} e^{-4ik^2(t-\tau)}   \left(2k g_0+ ig_1\right) \left(1+\frac{h_{21}}{k} + \frac{\hat{h}_{22}e^{-4ik^2 \tau}}{k^2}\right)    d\tau+O\left(\frac{1}{k^2}\right)+O\left(\frac{e^{-4ik^2 t}}
{k^2}\right) \\
 = & e^{-4ik^2 t}  \int_0^{t} -i \lambda \abs{g_0^2} \frac{\hat{h}_{11}(t)}{k} d\tau+ \int_0^{t} e^{-4ik^2(t-\tau)} \left( 2k g_0 \left(1+  \frac{\hat{h}_{22}e^{-4ik^2 \tau}}{k^2}\right)   \right)  d\tau\\
&+O\left(\frac{1}{k^2}\right)+O\left(\frac{e^{-4ik^2 t}}
{k^2}\right)\\
= & \frac{g_0(t)-g_0(0)e^{-4ik^2t}}{2ik} + e^{-4ik^2t}\int_0^{t}  \left( -i\lambda \abs{g_0^2} \frac{\hat{h}_{11}(t)}{k} + 2 g_0\frac{\hat{h}_{22}}{k}   \right)  d\tau+O\left(\frac{1}{k^2}\right)+O\left(\frac{e^{-4ik^2 t}}
{k^2}\right)\\
= & \frac{g_0(t)}{2ik}+\frac{ig_0(0)e^{-4ik^2t}}{2k}+O\left(\frac{1}{k^2}\right)+O\left(\frac{e^{-4ik^2 t}}
{k^2}\right).
\end{align*}
\end{itemize}
Comparing this with $O(1/k)+O(e^{-4ik^2 t}/k)$ contributions in (\ref{phi3}), it follows that for $\Phi_1(t,k)$,
\begin{align*}
h_{11}(t)      &=\frac{g_0(t)}{2i},\\
\hat{h}_{11}(t)&=\frac{ig_0(0)}{2},\\
h_{12}(t)      &=\frac{g_1(t)}{4}-\frac{i\lambda g_0(t)}{4}\int_{0}^{t} \overline{g_0(\tau)}g_1(\tau)-\overline{g_1(\tau)}g_0(\tau) d\tau,\\
\hat{h}_{12}(t)&=-\frac{g_1(0)}{4}-\frac{i\lambda g_0(0)}{4}\int_0^{t} \overline{g_0(\tau)}g_1(\tau)-\overline{g_1(\tau)}g_0(\tau) d\tau,\\
\end{align*}
and for $\Phi_2(t,k)$,
\begin{align*}
h_{21}(t)      &=\frac{\lambda}{2}\int_0^{t} \overline{g_0(\tau)}g_1(\tau)-\overline{g_1(\tau)}g_0(\tau) d\tau,\\
\hat{h}_{22}(t)&=\frac{\lambda\overline{g_0(t)}g_0(0)}{4i},\\
\hat{h}_{23}(t)&=\frac{i\lambda}{8}\left( \overline{g_1(t)}g_0(0)   -\overline{g_0(t)}g_1(0)-i\lambda g_0(0)\overline{g_0(t)}\int_0^{t} \overline{g_0(\tau)}g_1(\tau)-\overline{g_1(\tau)}g_0(\tau) d\tau \right).\\
\end{align*}
As a result, we have
\begin{align*}
A(k,T)=&\overline{\Phi_2(T,\overline{k})}\\
=& 1-\frac{\lambda}{2k}\int_{0}^{T}G(t)dt+\frac{i\lambda e^{4ik^2T}}{4k^2}\overline{g_0(0)}g_0(T)+\frac{-\frac{i\lambda}{8}\left( G(T)-i\lambda g_0(T)\overline{g_0(0)}\int_0^{T} G(\tau) d\tau \right)e^{4ik^2T}}{k^3}\\
&+O\left(\frac{1}{k^3}\right)+O\left(\frac{e^{-4ik^2 t}}
{k^4}\right),\\
B(k,T)=&-e^{4ik^2T}\Phi_1(T,k) \\
=& -\frac{g_0(T)e^{4ik^2 T}}{2ik}-\frac{ig_0(0)}{2k} -\frac{\left( g_1(T)-i\lambda g_0(T)\int_{0}^{T} G(\tau) d\tau \right) e^{4ik^2 T}}{4k^2}+\frac{g_1(0)+i\lambda g_0(0)\int_0^{T} G(\tau) d\tau}{4k^2}\\
&+O\left(\frac{1}{k^3}\right)+O\left(\frac{e^{-4ik^2 t}}
{k^3}\right),
\end{align*}
where
\[
G(t)=\overline{g_0(t)}g_1(t)-\overline{g_1(t)}g_0(t).
\]
\subsection*{Large $k$ asymptotics of $A(k,\infty),B(k,\infty)$}
We use the alternative set of equations for $A(k,\infty)=\tilde{\Phi}_2(0,k)$ and $B(k,\infty)=\tilde{\Phi}_1(0,k)$, which are defined using the linear Volterra integral equations (\ref{phi5}) and (\ref{phi6}).
Following similar steps to the calculations for $a(k)$ and $b(k)$, we have the expansions for $A(k,\infty),B(k,\infty)$,
\begin{align*}
A(k,\infty)&=\tilde{\Phi}_2(0,k)= 1-\frac{\lambda}{2k}\int_{0}^{\infty}G(t)dt+O\left(\frac{1}{k^2}\right),\\
B(k,\infty)&=\tilde{\Phi}_1(0,k)= \frac{g_0(0)}{2ik}+\frac{g_1(0)+i\lambda g_0(0)\int_0^{\infty} G(\tau) d\tau}{4k^2}+O\left(\frac{1}{k^3}\right),
\end{align*}
where
\[
G(t)=\overline{g_0(t)}g_1(t)-\overline{g_1(t)}g_0(t).
\]
These expansions are consistent with $A(k,T),B(k,T)$ by taking $T\rightarrow \infty$ after dropping all terms containing $e^{4ik^2 T}$. For $T=\infty$, the expansions are only valid in $\re{ik^2}\leq 0$.

\subsection*{Compatibility conditions and the expansions for the global relation}
If we expand the global relation at $k=\infty$, the compatibility conditions of the NLS equation at $x=0,t=0$ are obtained. In the case $T=\infty$, we get the global relation, $a(k)B(k,\infty)-b(k)A(k,\infty)=0$, at different orders:
\[O(1/k): \quad \frac{g_0(0)}{2ik}-\frac{q_0(0)}{2ik}=0, \quad
O(1/k^2): \quad \frac{g_1(0)}{4k^2}-\frac{q'_0(0)}{4k^2}=0.
\]
Although the expansions for $a,b,A,B$ depend on both $q_0(0),g_0(0),g_1(0)$ and the integrals of $q_0,g_0,g_1$, it turns out that the latter cancel in the expansion of the global relation. Furthermore, if $a,A\neq 0$, the global relation is rearranged as $b/a=B/A$, and
\[
B/A= \frac{g_0(0)}{2ik}+\frac{g_1(0)}{4k^2}+O\left(\frac{1}{k^3}\right),\quad  b/a= \frac{q_0(0)}{2ik}+\frac{q'_0(0)}{4k^2}+O\left(\frac{1}{k^3}\right).
\]
The integrals of $q_0,g_0,g_1$ do not appear in the expansions of $B/A$ and $b/a$.
For the UTM, the global relation is strictly satisfied throughout the calculation, but the solution does not require infinitely many compatibility conditions at the corner $x=t=0$. When the compatibility condition is violated at a certain order $1/k^n$, the unknown boundary function will be unbounded at the origin. This can be shown by an explicit example in the linear case.
\begin{example}
Using the UTM, the solution formula for the LS equation,
\[
iq_t+q_{xx}=0,
\]
on the half-line with the initial condition $q_0(x)=e^{-x}$ and the homogeneous Dirichlet boundary condition $g_0(t)=0$ is
\[
q(x,t)=\frac{1}{2\pi} \int_{-\infty}^{\infty} \frac{e^{ikx-ik^2t}}{1+ik}dk-\frac{1}{2\pi}\int_{\partial D^{+}}\frac{e^{ikx-ik^2t}}{1-ik}dk,
\]
where $\partial D^{+}$ is the boundary of the first quadrant, positively oriented. The Neumann data is to be computed
\[
g_1(t)=q_x(0,t)=\frac{1-i}{\sqrt{2\pi t}}-\frac{1}{\pi}\int_{-\infty}^{\infty}\frac{e^{-ik^2 t}}{1+k^2}dk.
\]
with asymptotics
\begin{align*}
q_x(0,t) = O(t^{-1/2}) &  \quad \text{ as  } t\rightarrow 0, \\
q_x(0,t) = O(t^{-3/2}) &  \quad \text{ as  } t\rightarrow \infty.
\end{align*}
Therefore, $g_1(t)=q_x(0,t)$ is unbounded at $t=0$. On the other hand, the global relation for $T=\infty$ is
\[
\hat{q}_0(k)=i\tilde{g}_1(k),   \quad \re{k}\leq 0, \im{k}\leq 0,
\]
with
\[
\hat{q}_0(k)=\int_0^{\infty} e^{-ikx} e^{-x}dx= \frac{1}{1+ik},
\]
\[
i\tilde{g}_1(k)=i\int_0^{\infty} e^{ik^2 t} q_x(0,t) dt=i\int_0^{\infty} \frac{ e^{ik^2 t} }{ \pi } \int_{-\infty}^{\infty} e^{-is^2 t} \frac{s^2}{1+s^2} ds  dt = \frac{1}{1+ik}.
\]
The global relation remains satisfied and the leading order expansion of $B(k,\infty)=i\tilde{g}_1(k)+k \tilde{g}_0(k)$ is given by
\[
B(k,\infty)=\lim_{t\rightarrow 0} \frac{\sqrt{2\pi}\left(q_x(0,t)\sqrt{t}\right)}{(i-1)k}+O(1/k^2).
\]
On the other hand, if $\lim_{t\rightarrow 0}\partial_x^n q(0,t)$ is known to be bounded, then the global relation implies that the initial and boundary values are compatible with respect to the LS equation up to the $n$-th order derivative.
\end{example}

\end{document}